\documentstyle[11pt]{article}
\input amssym.def

\newtheorem{df}{Definition}[section]
\newtheorem{thm}[df]{Theorem}
\newtheorem{prop}[df]{Proposition}
\newtheorem{cor}[df]{Corollary}
\newtheorem{ex}[df]{Example}
\newtheorem{lem}[df]{Lemma}

\def\boom{\quad\lower3pt\hbox{\vrule height1.1ex width .9ex depth -.2ex}
                    \vskip9pt}
\def\kwo#1{\lq\lq #1\rq\rq}
\def\st{\ \vert\ }
\newcommand{\pf}{\noindent{\sc Proof.}\ }

\let\Bar=\overline
\let\Hat=\widehat
\let\Tilde=\widetilde
\let\Vec=\overrightarrow

\def\tilalpha{\widetilde\alpha}
\def\tilbeta{\skew6\widetilde\beta}
\def\tiliota{\widetilde\iota}
\def\til0{\widetilde 0}
\def\tilone{\widetilde 1}
\def\tila{\widetilde a}
\def\tilq{\widetilde q}

\let\da=\partial
\let\isom=\cong

\def\upa{\uparrow}

\def\cls{\clubsuit}
\def\hs{\heartsuit}
\def\ds{\diamondsuit}
\def\spades{\spadesuit}

\def\chigh{{\raise1.5pt\hbox{$\chi$}}}
\let\Ga=\Gamma
\let\Om=\Omega
\let\Th=\Theta

\def\LAgpd{${\cal LA}$-groupoid}
\def\LAdgpd{${\cal LA}$-double groupoid}

\def\VBgpd{${\cal VB}$-groupoid}

\def\CDO{\mathop{\rm CDO}}
\def\Exp{\mathop{\rm Exp}}
\def\id{{\rm id}}

\def\BarH#1{\Bar{#1}^{\scriptscriptstyle H}}
\def\BarV#1{\Bar{#1}^{\scriptscriptstyle V}}
\def\VecH#1{\Vec{#1}{\rule{0cm}{1ex}}^{\scriptscriptstyle H}}

\def\ddt#1{\left.\frac{d}{dt}#1\right|_0}
\def\proj#1{\stackrel{#1}{\sim}}
\def\riri{\Vec{\phantom{X}}}
\def\ts#1{\stackrel{\vee}{#1}}
\def\twiddle{\!\Tilde{\phantom{X}}}

\def\inj{>\!\!\!-\!\!\!-\!\!\!-\!\!\!>}
\def\surj{-\!\!\!-\!\!\!-\!\!\!\gg}
\def\pback#1{\mathbin{{{\lower1.2ex\hbox{$\times$}}\atop #1}}}
\def\sdp{\mathbin{\hbox{$\mapstochar\kern-.3333em\times$}}}

\def\gpd{\,\lower1pt\hbox{$\longrightarrow$}\hskip-.24in\raise2pt
             \hbox{$\longrightarrow$}\,}

\def\sgpd{\,\lower1pt\hbox{$\mlra$}\hskip-0.4in\raise2pt\hbox{$\mlra$}\,}
\def\svgpd{\Big\downarrow\!\!\Big\downarrow}

\def\vgpd{\Bigg\downarrow\!\!\Bigg\downarrow}

\def\vlra{\hbox{$\,-\!\!\!-\!\!\!-\!\!\!-\!\!\!-\!\!\!
-\!\!\!-\!\!\!-\!\!\!-\!\!\!-\!\!\!\longrightarrow\,$}}

\def\lgpd{\,\lower1pt\hbox{$\vlra$}\hskip-1.02in\raise2pt\hbox{$\vlra$}\,}

\def\lrah{\hbox{$\,-\!\!\!-\!\!\!
-\!\!\!-\!\!\!-\!\!\!-\!\!\!-\!\!\!\longrightarrow\,$}}

\def\mlra{\hbox{$\,-\!\!\!-\!\!\!\longrightarrow\,$}}

\def\hcompo#1#2#3{{\vcenter{\vbox{\hrule height.#2pt\hbox{\vrule width.#2pt
   height#1pt\kern#3pt\vrule width.#2pt\kern#3pt\vrule width.#2pt}
   \hrule height.#2pt}}}}

\def\hcomp{\mathchoice\hcompo643\hcompo643\hcompo{2.2}3{1.1}
                       \hcompo{1.4}3{.7}}

\def\vcompo#1#2#3{{\vcenter{\vbox{\hrule height.#2pt
   \hbox{\vrule width.#2pt height#3pt\kern#1pt\vrule width.#2pt}
   \hrule height.#2pt
   \hbox{\vrule width.#2pt height#3pt\kern#1pt\vrule width.#2pt}
   \hrule height.#2pt}}}}

\def\vcomp{\mathchoice\vcompo643\vcompo643\vcompo{2.2}3{1.1}
                       \vcompo{1.4}3{.7}}

\def\dsq{\mathop{\lower1pt\vbox{\hrule height.4pt \hbox
{\vrule width.4pt height.6em
\kern.6em \vrule width.4pt} \hrule height.4pt}}}

\def\dcomp{\mathop{\dsq\hskip-.88em\raise1pt\hbox{$\scriptstyle\nwarrow$}}}

\def\ssq{\vbox{\hrule height.4pt \hbox{\vrule width.4pt height.7in
\kern.7in \vrule width.4pt} \hrule height.4pt}}

\def\tsq{\mathop{\lower1pt\vbox{\hrule height.4pt \hbox
{\vrule width.4pt height.7em
\kern.7em \vrule width.4pt} \hrule height.4pt}}}

\def\dpl{+\hskip-6pt +\hskip4pt}

\def\dtimes{\mathbin{\hbox{\huge.}}}       

\def\affel#1#2#3#4#5#6#7#8{\matrix{#7&\raise.7ex\hbox{$#1$}&#6\cr
                  \raise.5in\hbox{$#2$}&\sq&\raise.5in\hbox{$#3$}\cr
                          #8&#4&#5\cr}}

\def\dvb#1#2#3#4#5#6#7#8{\matrix{&&#5&&\cr
                                 &#1&\lrah&#3&\cr
                                 &&&&\cr
                               #6&\Big\downarrow&&\Big\downarrow&#7\cr
                                 &&&&\cr
                                 &#2&\lrah&#4&\cr
                                 &&#8&&\cr}}

\def\R{\Bbb R}

\begin{document}

\title{{\bf DOUBLE LIE ALGEBROIDS AND SECOND-ORDER GEOMETRY, II}
}

\author{K. C. H. Mackenzie\\
        School of Mathematics and Statistics\\
        University of Sheffield\\
        Sheffield, S3 7RH\\
        England\\
        {\sf K.Mackenzie@sheffield.ac.uk}}

\date{{\sf April 26, 1997}}

\maketitle

\begin{abstract}
We complete the construction of the double Lie algebroid of a double
Lie groupoid begun in the first paper of this title. We extend the
construction of the tangent prolongation of an abstract Lie algebroid to
show that the Lie algebroid structure of any \LAgpd\ may be prolonged
to the Lie algebroid of its groupoid structure. In the case of a double
groupoid, this prolonged structure for either \LAgpd\ is canonically
isomorphic to the Lie algebroid structure associated with the other; this
extends many canonical isomorphisms associated with iterated tangent and
cotangent structures.

We calculate several examples from Poisson geometry. We show that the
cotangent of any double Lie groupoid is a symplectic double groupoid and
that the side groupoids of a symplectic double groupoid are Poisson groupoids
in duality; thus the duals of the \LAgpd s of any double groupoid are a pair
of Poisson groupoids in duality.
\end{abstract}

\newpage

In the first paper of this title \cite{Mackenzie:1992} we proposed a Lie
theory for double Lie groupoids based on the known Lie theory of
ordinary Lie groupoids, on various constructions in Poisson geometry, and
on the special features of the differential geometry of the tangent bundle.
In this paper we complete the lengthy construction of the double Lie algebroid
of a double Lie groupoid which was begun in \cite{Mackenzie:1992}, and
calculate the fundamental examples, concentrating on those which arise in
Poisson geometry. We refer to the introduction of \cite{Mackenzie:1992} for
an overview of the programme and its ultimate aims. Here we describe the
features specific to the present paper, by way of orientation.

Ordinary Lie algebroids may be viewed both as generalizations of Lie
algebras---and therefore as vehicles for a generalized Lie theory
\cite{Mackenzie:1995}---and as abstractions of the tangent bundle of an
ordinary manifold. (For the latter point of view see, for example,
\cite{Weinstein:1996} or \cite{MackenzieX:1998}.) For double Lie algebroids,
the model of the Lie theory of Lie groups and Lie algebras is a more distant
one, and we will be chiefly concerned here with constructions which derive
ultimately from the calculus possible on the double tangent bundle.

In \cite{Mackenzie:1992} it was shown how a single application of the Lie
functor to a double Lie groupoid $S$ produces an \LAgpd, that is, a
Lie groupoid object in the category of Lie algebroids. If one applies the
Lie functor to, say, the vertical structure of $S$, then the \LAgpd\ is
vertically a Lie algebroid and horizontally a Lie groupoid; this follows
easily from the fact that the Lie functor preserves pullbacks and the
diagrams which define a groupoid structure. One may then take the Lie
algebroid of this horizontal groupoid and obtain a double vector bundle
whose horizontal structure is a Lie algebroid. Interchanging the order of
the processes yields a second double vector bundle with the Lie algebroid
structure now placed vertically. The two double vector bundles may be
identified by a map derived from the canonical involution in the double
tangent bundle of $S$, and one thus obtains a double vector bundle all four
sides of which have Lie algebroid structures; this is the double Lie
algebroid of $S$. This is the construction as outlined in the introduction
to \cite{Mackenzie:1992}. It includes, amongst many others, the double
tangent bundle of an arbitrary manifold, the Lie bialgebra of a Poisson
group, and the double cotangent bundle of a Poisson manifold.

However, abstract \LAgpd s frequently arise in nature without an underlying
double Lie groupoid, most notably the cotangent \LAgpd\ of a Poisson groupoid
or Poisson Lie group. We have therefore shown in \S\ref{sect:dlalag} here
that the Lie algebroid structure on an \LAgpd\ may be prolonged to the Lie
algebroid of the groupoid structure. This extends, in particular, the
construction of tangent Lie algebroid structures given by Ping Xu and the
author in \cite{MackenzieX:1994}. We must then prove (Theorem \ref{thm:jS2})
that, in the case of the \LAgpd s of a double Lie groupoid, the
prolonged Lie algebroid structure on the double Lie algebroid of either
\LAgpd\ coincides with the Lie algebroid of the Lie groupoid structure
of the other. This result embodies many canonical isomorphisms known in
special cases.

We believe that the basic simplicity of this construction, and the richness
of the phenomena which it encompasses, establish beyond doubt that it is a
natural and correct definition. On a future occasion we will consider the
abstract notion of double Lie algebroid and its local integrability.
In \S\ref{sect:aff} to \S\ref{sect:sdgapdg} of the present paper, we are
concerned to demonstrate the range of examples covered by the notion of
double Lie algebroid.

The paper concludes with a detailed examination of the example of
symplectic and Poisson double groupoids, and the duality associated with
the former. We give in \S\ref{sect:datcdg} a double cotangent groupoid
construction analogous to the cotangent of an ordinary groupoid
\cite{CDW}, and show that this embodies the relationship between the
groupoid structures and the symplectic structure. The key to the cotangent
double groupoid is the notion of core introduced in \cite{BrownM:1992}.
In \S\ref{sect:sdgapdg} we give a simple proof that the side groupoids of
a symplectic double groupoid are Poisson groupoids the Lie bialgebroids of
which are dual, and that the core gives a symplectic realization of the
double base. The double Lie algebroid is (for either structure) the
cotangent Lie algebroid of the linearized Poisson structures associated to
the dual Poisson groupoids; this is the object studied at length in
\cite{MackenzieX:1994} and elsewhere as embodying the properties of a Lie
bialgebroid.

In \S\ref{sect:aff} we calculate the double Lie algebroid associated with
an affinoid structure: this is, as indicated by Weinstein
\cite{Weinstein:1990}, a pair of conjugate flat partial connections.
We show that the infinitesimal version of the notion of butterfly diagram is
equivalent (under mild conditions) to such a pair of connections.

We have kept to a minimum the repetition of material from
\cite{Mackenzie:1992} and so refer the reader there for basic definitions
and notations and much else. We also make extensive use of the material in
\S5 and \S7 of \cite{MackenzieX:1994}, and of \cite{MackenzieX:1998}. One
minor change from \cite{Mackenzie:1992} is that the pullback of a vector
bundle $A$ over a map $f$ is denoted $f^!A$ rather than $f^*A$,
with other pullback notation modified accordingly.

An early version of much of the material of this paper was presented in May
1993 at a session of the Paris S\'eminaire de G\'eom\'etrie et M\'ecanique
under the direction of Professor C.-M. Marle. I am grateful to him for that
opportunity, and to Johannnes Huebschmann for organizing the context in which
that visit took place. Most of the detailed writing was done amid the
admirable facilities and working conditions of the Centre \'Emile Borel in
a stay of five weeks during the Semestre de G\'eom\'etrie Symplectique, and
I acknowledge with gratitude the support of the Royal Society and of the
Ecole Polytechnique which made that visit possible.

\section{\bf THE DOUBLE LIE ALGEBROID OF AN ${\cal LA}$-GROUPOID}
\label{sect:dlalag}

The construction of the double Lie algebroid of a double Lie groupoid will
be made much easier if we begin by first considering the application of the
Lie functor to an arbitrary \LAgpd
\begin{figure}[htb]
\begin{picture}(395,100)
\put(130,50){$\matrix{&&\tilalpha,\tilbeta &&\cr
        &\Omega  &\lgpd&A &\cr
        &&&&\cr
  \tilq &\Bigg\downarrow&&\Bigg\downarrow&\cr
        &&&&\cr
        &G  &\lgpd&M,&\cr
        &&\alpha,\beta&&\cr}$}
\put(180,50){$TG\stackrel{T(\alpha),T(\beta)}{\lgpd}\ TM$}
\put(155,75){\vector(1,-1){15}}    \put(170,65){$\tila$}
\put(265,75){\vector(1,-1){15}}
\put(175,45){\vector(-1,-1){15}}
\put(285,45){\vector(-1,-1){15}}
\end{picture}         \caption{\ \label{fig:LAgpd}}
\end{figure}
as defined in \cite[\S4]{Mackenzie:1992} (see Figure~\ref{fig:LAgpd}). Here
$\Om$ is both a Lie
groupoid over base $A$ (which is itself a Lie algebroid over $M$), and a Lie
algebroid over $G$ (which is a Lie groupoid over $M$); the two structures
on $\Om$ commute in the sense that the maps defining the groupoid structure
are all Lie algebroid morphisms. It is further assumed that the double source
map $(\tilq,\tilalpha)\colon\Om\to G\pback{M}A$ is a surjective submersion.

Since $\Om$ is a Lie groupoid over $A$ we may take its Lie algebroid $A\Om$;
denote the bundle projecttion $A\Om\to A$ by $\ts{q}$.
Since $\tilq\colon\Om\to G$ is a morphism of groupoids, it induces a morphism
of Lie algebroids $A(\tilq)\colon A\Om\to AG$. And since each of the other
maps defining the vector bundle structure in $\Om\to G$ (namely the addition,
scalar multiplication, and zero section) are morphisms of Lie groupoids,
we may apply the Lie functor to them, and we obtain a vector bundle structure
on $A(\tilq)\colon A\Om\to AG$. As with the construction in
\cite[\S4, p.216]{Mackenzie:1992},
of which this is a special case, a detailed proof is unnecessary: a vector
bundle is defined by pullback diagrams, and the Lie functor preserves
pullbacks and diagrams.

\begin{prop}
The above construction yields a double vector bundle
\begin{equation}                               \label{eq:AOm}
\matrix{&&\ts{q}&&\cr
        &A\Om&\lrah &A&\cr
        &&&&\cr
A(\tilq)&\Bigg\downarrow& &\Bigg\downarrow&q\cr
        &&&&\cr
        &AG&\lrah &M.&\cr
        &&q_G&&\cr}
\end{equation}
\end{prop}

\pf
To verify the commutativity, it only needs to be recalled that a vector
bundle, regarded as a Lie groupoid, is itself its own Lie algebroid, and,
similarly, applying the Lie functor to a morphism of vector bundles leaves the
morphism unchanged.
\boom

Three sides of (\ref{eq:AOm}) have Lie algebroid structures, and we are now
concerned to put such a structure on $A\Om\to AG$ as well. This will take
us until the end of the section. The method used is a generalization of the
construction of a tangent Lie algebroid \cite[5.1]{MackenzieX:1994}.

The anchor is the composite of $A(\tila)\colon A\Om\to ATG$ with the inverse
of the canonical isomorphism $j_G\colon TAG\to ATG$ of \cite[7.1]{MackenzieX:1994}; we
denote it by ${\bf a}$. Using \cite[7.1]{MackenzieX:1994}, ${\bf a}$ is a morphism of Lie
algebroids over $a\colon A\to TM$.

We define the bracket structure on $A\Om\to AG$ in terms of sections of two
specific types.

\begin{df}
A {\em star section} of the \LAgpd\ in Figure~\ref{fig:LAgpd} is a pair
$(\xi,X)$ where $\xi\in\Ga_G\Om$, $X\in\Ga A$, and we have $\tilone\circ X =
\xi\circ 1,\ \tilalpha\circ\xi = X\circ\alpha.$
\end{df}

This is a weakened form of the notion of morphic section considered in
\cite[4.3]{Mackenzie:1992}.
The terminology comes from \cite[3.2]{HigginsM:1990a}.

\begin{lem}                            \label{lem:ss}
Given any $X\in\Ga A$, there is a star section $(\xi,X)$ of the \LAgpd\
in Figure~\ref{fig:LAgpd}.
\end{lem}

\pf
Define $\eta\in\Ga_G\Om$ by setting $\eta(1_m) = \tilone_{X(m)}$ for
$m\in M$, and extending over $G$. Then $\mu\colon G\to A$ defined by
$\mu(g) = \tilalpha(\eta(g)) - X(\alpha g)$ is a section of the pullback
bundle $\alpha^!A$. Thanks to the double source condition, there is a
$\zeta\in\Ga_G\Om$ with $\tilalpha\circ\zeta = \mu$; we can also require that
$\zeta$ vanish on all $1_m\in G$. Now $\xi = \eta - \zeta$ is a section
of the required type.
\boom

The conditions on a star section ensure that it is possible to apply
the Lie functor to it; we obtain a linear section $A(\xi)\colon AG\to A\Om$
of $A(\tilq)$ for which $\ts{q}\circ A(\xi) = X\circ q_G$. As in
\cite[4.3]{Mackenzie:1992} it is
clear that the bracket of two star sections is a star section. We now define
\begin{equation}                                      \label{eq:bracket1}
[A(\xi),A(\eta)] = A([\xi,\eta])
\end{equation}
for any two star sections $(\xi,X),\ (\eta,Y)$.

To define sections of the second type, we must make a detour to determine the
core structure of the double
vector bundle (\ref{eq:AOm}) in terms of the core of the underlying \LAgpd.
The core Lie algebroid of an \LAgpd\ was defined in
\cite[\S5]{Mackenzie:1992}.

\begin{prop}
Let $(\Om;G,A;M)$ be an \LAgpd\ with core Lie algebroid $K$.
Denote by $K_{vb}$ the vector bundle underlying $K$. Considering $K_{vb}$
as a groupoid, let it act on $q_A\colon A\to M$ by
$(\kappa,X)\mapsto \da_A(\kappa) + X,$ where
$\kappa\in K, X\in A,\ q_K(\kappa) = q_A(X)$. Then
\begin{equation}                                    \label{eq:LAgpdes}
K_{vb}\sdp q_A \inj \Om \stackrel{\tilq}{\surj} G
\end{equation}
is an exact sequence of Lie groupoids, where the injection is
$(\kappa,X)\mapsto \kappa + \tilone_X$.

Let the Lie algebroid $K$ act on $\beta\colon G\to M$ by sending $\kappa\in
\Ga K$ to $\Vec{\tila(\kappa)}\in\Ga TG$, the right-invariant vector
field corresponding to $\tila(\kappa)\in\Ga AG$. Then
\begin{equation}                              \label{eq:LAgpdcore}
K\sdp\beta \inj \Om \stackrel{\tilalpha}{\surj} A
\end{equation}
is an exact sequence of Lie algebroids, where the injection is $(\kappa,g)
\mapsto \kappa\til0_g$.
\end{prop}

\pf
Take $(\kappa_1,X_1), (\kappa_2,X_2)\in K\sdp q_A$ with
$X_1 = \da_A(\kappa_2) + X_2$. Multiplying the images in $\Om$, and using the
interchange law between the vector bundle and groupoid structures, we have
\begin{eqnarray*}
(\kappa_1 + \tilone_{X_1})(\kappa_2 + \tilone_{X_2}) & = &
(\kappa_1 + \tilone_{X_2} +
\tilone_{\da_A(\kappa_2)})(\tilone_{X_2} + \kappa_2)\\
& = &
(\kappa_1 + \tilone_{X_2})\tilone_{X_2} + \tilone_{\da_A(\kappa_2)}\kappa_2\\
&  = & \kappa_1 + \tilone_{X_2} + \kappa_2,
\end{eqnarray*}
which is the image of $(\kappa_1 + \kappa_2, X_2)$ under the injection. So
the injection is a groupoid morphism. Exactness is easily verified.

The injection in the second statement preserves the anchors since
$$
\tila(\kappa\til0^\Om_g) = \tila(\kappa)\til0^{TG}_g = T(R_g)(\tila(\kappa))
= \Vec{\tila(\kappa)}(g),
$$
using the fact that $\tila$ is a morphism of groupoids. It suffices to
verify bracket preservation on sections of the form $1\otimes\kappa\in
\Ga(K\sdp\beta)$, where $\kappa\in\Ga K$. This follows from (22) in
\cite[\S5]{Mackenzie:1992}.
\boom

For the \LAgpd\ $(TG;G,TM;M)$ of a Lie groupoid $G$, these two sequences are
$$
A_{vb}G\sdp p_M \inj TG \stackrel{p_G}{\surj} G \quad\mbox{and}\quad
AG\sdp\beta \stackrel{{\cal R}}{\inj} TG \stackrel{T(\alpha)}{\surj} TM,
$$
where ${\cal R}$ is the right translation map $(X,g)\mapsto T(R_g)(X)$.

For completeness, we state the corresponding results for double Lie
groupoids.

\begin{prop}
Let $(S;H,V;M)$ be a double Lie groupoid with core groupoid $C$ (see
\cite{BrownM:1992} or \cite[\S2]{Mackenzie:1992}). Then
\begin{equation}             \label{eq:coreseqV}
C\sdp\beta_H\inj S_V \stackrel{\tilalpha_H}{\surj} V
\end{equation}
is an exact sequence of Lie groupoids, where $C$ acts on $\beta_H$ by
$(k,h) \mapsto \da_H(k)h$ and the injection is
$(k,h)\mapsto k\,\hcomp\,\tilone^V_h$.

Likewise
\begin{equation}             \label{eq:coreseqH}
C\sdp\beta_V\inj S_H \stackrel{\tilalpha_V}{\surj} H
\end{equation}
is an exact sequence of Lie groupoids, where $C$ acts on $\beta_V$ by
$(k,v) \mapsto \da_V(k)v$ and the injection is
$(k,v)\mapsto k\,\vcomp\,\tilone^H_v$.
\end{prop}

When $S$ is a locally trivial double Lie groupoid, the exact sequences
(7), (8) of \cite[\S2]{Mackenzie:1992} follow easily from
(\ref{eq:coreseqV}), (\ref{eq:coreseqH}).

\begin{prop}                               \label{prop:cores}
Let $(S;H,V;M)$ be a double Lie groupoid with core groupoid $C$. Then
$AC$, the Lie algebroid of $C$, identifies canonically with the core
Lie algebroid of the \LAgpd\ $A_HS$ and with the core Lie algebroid of
the \LAgpd\ $A_VS$.
\end{prop}

\pf
Let $K$ be the core of $A_HS$. Fix $m\in M$. Then a curve $c_t$ in the
$\alpha$-fibre of $C$ above $m$ with $c_0 = 1^C_m$ is also a curve in the
$\alpha$-fibre of $S_H$ above $1^V_m$, with $c_0 = 1^2_m$, the double
identity at $m$. Thus there is a vector bundle morphism $AC\to A_HS$, which
is easily seen to be an isomorphism over $M$ onto $K$. Now applying the Lie
functor to (\ref{eq:coreseqH}) and comparing it with (\ref{eq:LAgpdcore})
shows that this is an isomorphism of Lie algebroids.
\boom

Returning to the \LAgpd\ in Figure~\ref{fig:LAgpd}, we can now apply the
Lie functor to (\ref{eq:LAgpdes}), and obtain an exact sequence of Lie
algebroids
\begin{equation}                                    \label{eq:LAgpdesLA}
K^0\sdp q_A \inj A\Om \stackrel{A(\tilq)}{\surj} AG.
\end{equation}
Here $K^0$ is the vector bundle underlying $K$ equipped with the zero
anchor and zero bracket: it is the Lie algebroid of $K_{vb}$. The action
of $K^0$ on $q_A\colon A\to M$ is $\kappa^\dagger = \da_A(\kappa)^\upa$,
where for any $X\in\Ga A$ we denote by $X^\upa$ the vertical vector field
on $A$ defined by
$$
X^\upa(F)(Y) = \ddt{F(Y+tX(qY))}
$$
where $F\in C(A)$ and $Y\in A$. In \cite[(23)]{MackenzieX:1994}, $X^\upa$
is denoted $\breve{X}$ and called the {\em core vector field}, but here we
use the notation $X^\upa$; in the case $A = TM$, $X^\upa$ is the vertical
lift of a vector field on $M$ to a vector field on $TM$.

If we disregard the Lie algebroid structures in (\ref{eq:LAgpdesLA}), we
can rewrite it as an exact sequence of vector bundles over $A$:
\begin{equation}                                    \label{eq:LAgpdesLAvb}
q_A^!K \inj A\Om \stackrel{A(\tilq)^!}{\surj} q_A^!AG.
\end{equation}
From general properties of vector bundles
\cite[Proposition 1.2]{Mackenzie:1992}, the next result then follows.

\begin{prop}                           \label{prop:dvb}
The core of the double vector bundle (\ref{eq:AOm}) is $K_{vb}$, the vector
bundle underlying the core of $\Om$.
\end{prop}

We think of elements of $A\Om$ as derivatives
$$
\Xi = \ddt{\xi_t}
$$
where $\xi_t$ is a curve in a fixed $\alpha$-fibre $\Om_X$ of $\Om$ with
$\xi_0 = \tilone_X$. Thus $\ts{q}\!(\Xi) = X$. If we write
$X_t = \tilbeta(\xi_t)$ and $g_t = \tilq(\xi_t)$ then
$$
\ddt{X_t} = \ \ts{a}\!(\Xi)\in TA\qquad\mbox{and}\qquad
\ddt{g_t} = A(\tilq)(\Xi)\in AG.
$$
If $k$ is an element of the core $K$ of $\Om$, then it has
$\tilalpha(k) = 0_m,\ \tilq(k) = 1_m$ and $\tilbeta(k) = Y$ for some
$m\in M$ and $Y = \da_A(k)\in A$. The corresponding element of $A\Om$
is $\ddt{(tk)}$, where $tk$ is the scalar multiplication in the bundle
$\Om\to G$. We denote this by $\Bar{k}$ if confusion is likely. Note that,
using the notation of \cite[\S5]{MackenzieX:1994}, $\ts{a}\!(\Bar{k}) = \Bar{Y}$, the
core element of $TA$ corresponding to $Y\in A$.

The injection in (\ref{eq:LAgpdesLA}) can now be written in either of two
ways:
$$
(k,X)\mapsto \ddt{(tk + \tilone_X)} = \Bar{k}\dpl \til0_X
$$
where $\dpl$ denotes the addition in $A\Om\to AG$ and $\til0_X$ is the
zero in $A\Om\to A$ above $X$. Given a section $\kappa\in\Ga K$, we define
a section $\kappa^\ds$ of $A\Om\to A$ by
\begin{equation}                               \label{eq:coreds}
\kappa^\ds(X) = \Bar{\kappa(q_AX)}\dpl \til0_X = \ddt{(\tilone_X +
t\kappa(q_AX))}.
\end{equation}
We call $\kappa^\ds$ the {\em core section} of $A\Om\to A$ corresponding to
$\kappa$.

\begin{prop}                        \label{prop:upa}
For $\kappa, \lambda\in\Ga K$ and $f\in C(M)$ we have
$$
(\kappa+\lambda)^\ds = \kappa^\ds + \lambda^\ds,\qquad
(f\kappa)^\ds = (f\circ q_A)\kappa^\ds,\qquad
[\kappa^\ds,\lambda^\ds] = 0,\qquad
\ts{a}(\kappa^\ds) = \da_A(\kappa)^\upa.
$$
\end{prop}

\pf
The first two statements are trivial. For the third, note that $\kappa^\ds$
is the image under the injection in (\ref{eq:LAgpdesLA}) of the pullback
section $1\otimes\kappa$ of $q_A^!K$. In any action Lie algebroid one has
$[1\otimes\kappa,1\otimes\lambda] = 1\otimes[\kappa,\lambda]$
\cite[2.4]{HigginsM:1990a}, and since the acting Lie algebroid is $K^0$,
we have $[\kappa,\lambda] = 0$.

Since the injection in (\ref{eq:LAgpdesLA}) is a Lie algebroid morphism,
$\ts{a}(\kappa^\ds(X)) = a'(\kappa(m),X)$, where $a'$ is the anchor
for $K^0\sdp q_A$ and $m = q_AX$. But we know that the anchor of an action
Lie algebroid is given by the action, and so $a'(\kappa(m),X) =
\kappa^\dagger(X) = \da_A(\kappa)^\upa(X)$.
\boom

\begin{prop}                 \label{prop:anchors}
\begin{enumerate}
\item
The anchor $\ts{a}\colon A\Om\to TA$ is a morphism of double vector bundles
over $a_G\colon AG\to TM$ and $\id_A$ with core morphism
$\da_A\colon K\to A$.

\item
The map ${\bf a}\colon A\Om\to TAG$ is a morphism of double vector bundles
over $a\colon A\to TM$ and $\id_{AG}$ with core morphism
$\da_{AG}\colon K\to AG$.
\end{enumerate}
\end{prop}

\pf
Only the statements concerning the core morphisms require proof. For the
first statement the proof is essentially the same as for the last
result of Proposition~\ref{prop:upa}. For the second, it suffices to prove
that $A(\tila)\colon A\Om\to ATG$ has core morphism $\da_{AG}$. Now
for $k\in K$,
$$
A(\tila)\left(\ddt{(tk)}\right) = \ddt{t\tila(k)} = \ddt{t\da_{AG}(k)} =
\Bar{\da_{AG}(k)}.
$$
\boom

From Proposition \ref{prop:dvb} it follows that there is an exact
sequence of vector bundles over $AG$, a companion to (\ref{eq:LAgpdesLA}),
\begin{equation}                                    \label{eq:LAgpdesLAvb2}
q_G^!K \inj A\Om \stackrel{\ts{q}^!}{\surj} q_G^!A,
\end{equation}
where the injection is $(k,x)\mapsto A(0)(x) + \Bar{k}$.

Again, given a section $\kappa$ of $K$, there is a section $\kappa^\hs$
of $A\Om\to AG$ defined by
\begin{equation}                           \label{eq:coresection}
\kappa^\hs(x) = A(0)(x) + \Bar{\kappa(q_{AG}(x))}
\end{equation}
for $x\in AG$. We call $\kappa^\hs$ the {\em core section} of $A\Om\to
AG$ corresponding to $\kappa$. The construction $\kappa\mapsto\kappa^\hs$
generalizes the construction $X\mapsto \Hat{X}$ for a vector bundle in
\cite[\S5]{MackenzieX:1994}.

\begin{prop}
For $\kappa, \lambda\in\Ga K$ and $f\in C(M)$ we have
$$
(\kappa+\lambda)^\hs = \kappa^\hs + \lambda^\hs,\qquad
(f\kappa)^\hs = (f\circ q_{AG})\kappa^\hs,\qquad
{\bf a}\circ\kappa^\hs = \da_{AG}(\kappa)^\upa.
$$
\end{prop}

\pf The first two statements are trivial, and the last follows from
Proposition \ref{prop:anchors}(ii).
\boom

It is the core sections $\kappa^\hs,\ \kappa\in\Ga K$, that together with
those of the form $A(\xi)$ for $(\xi,x)$ a star section, will allow us to
define the Lie algebroid structure on $A\Om\to AG$.

\begin{prop}                          \label{prop:generation}
The sections of the form $A(\xi)$, where $(\xi,X)$ is a star section, and
$\kappa^\hs$, where $\kappa\in\Ga K$, generate $\Ga_{AG}A\Om$.
\end{prop}

\pf
Take $\Xi\in A\Om$ with $A(\tilq)(\Xi) = x$ and $\ts{q}\!(\Xi) = X(m)$. Extend
$X(m)\in A$ to a section $X$ of $A$. By Lemma \ref{lem:ss}, there is a star
section $(\xi,X)$ of $\Om$. We now have $\ts{q}\!(A(\xi)(x)) = X(m)$ and
so, by (\ref{eq:LAgpdesLAvb2}), we have
$$
\Xi = A(\xi)(x) \dpl \kappa^\hs(x)
$$
for some $\kappa\in\Ga K$.
\boom

For $\kappa,\lambda \in\Ga K$ we define
\begin{equation}                             \label{eq:bracket2}
[\kappa^\hs,\lambda^\hs] = 0.
\end{equation}
It remains to define brackets of the form $[A(\xi),\kappa^\hs]$. We do
this by an extension of the methods developed in \cite{MackenzieX:1998}.

For $\kappa\in\Ga K$, denote (as in \cite[\S5]{Mackenzie:1992}) by
$\Bar{\kappa}$ the element of
$\Ga_G\Om$ given by $\Bar{\kappa}(g) = \kappa(\beta g)\til0_g$. Then
$\Bar{\kappa}\stackrel{\alpha}{\sim}0$. If $(\xi,X)$ is a star section, we
have $\xi\stackrel{\alpha}{\sim}X$ and so it follows that
$[\xi,\Bar{\kappa}]\stackrel{\alpha}{\sim}0$. Define
\begin{equation}                                 \label{eq:D}
D_\xi(\kappa) = [\xi,\Bar{\kappa}]\circ 1\in\Ga K.
\end{equation}

\begin{lem}                           \label{lem:D}
Let $(\xi,X)$ and $(\eta,Y)$ be star sections of $\Om$ and let $\kappa$
be a section of $K$. Then:
\begin{enumerate}
\item If $\zeta\in\Ga_G\Om$ is any section with $\zeta\circ 1 = \kappa$,
we have $D_\xi(\kappa) = [\xi,\zeta]\circ 1.$

\item $D_{[\xi,\eta]}(\kappa) = D_\xi(D_\eta(\kappa)) -
D_\eta(D_\xi(\kappa)).$

\item For any $f\in C(M)$, $D_\xi(f\kappa) = fD_\xi(\kappa) + a(X)(f)\kappa.$
\end{enumerate}
\end{lem}

\pf
(i) It only needs to be proved that if $\eta\in\Ga_G\Om$ has
$\eta\circ 1 = 0$ then $[\xi,\eta]\circ 1 = 0$. But $\tilone\colon A\to\Om$
is a Lie algebroid morphism over $1\colon M\to G$, and
$X\stackrel{1}{\sim}\xi$, so it follows that if $0\stackrel{1}{\sim}\eta$,
then $0 = [X,0]\stackrel{1}{\sim}[\xi,\eta],$ which is to say that
$[\xi,\eta]\circ 1 = 0.$ From this (ii) now follows immediately.

(iii) For each $m\in M$ we have
\begin{eqnarray*}
D_\xi(f\kappa)(m) & = & [\xi,(f\circ\beta)\Bar{\kappa}](1_m)\\
& = & f(m)D_\xi(\kappa)(m) + \tila(\xi)(1_m)(f\circ\beta)\kappa(m)
\end{eqnarray*}
and $\tila(\xi(1_m)) = \tila(\tilone_{X(m)}) = T(1)(a(X(m)))$ since
$\tila\colon\Om\to TG$ is a groupoid morphism over $a\colon A\to TM$.
\boom

Thus $D_\xi\colon\Ga K\to\Ga K$ is a covariant differential operator over
$a(X)$.

\begin{prop}
Let $(\xi,X)$ be a morphic section. Then $[\xi,\Bar{\kappa}] =
\Bar{D_\xi(\kappa)}$ for all $\kappa\in\Ga K$, and $D_\xi$ is a derivation
of the bracket structure on $\Ga K$.
\end{prop}

\pf
Recall from \cite[\S5, p.230]{Mackenzie:1992}, that a section $\zeta$ of
$\Bar{K}$, the kernel of $\tilalpha\colon\Om\to A$, is $G$-equivariant
if $\zeta(gh) = \zeta(g)\til0_h$
for all compatible $g,h$. A $G$-equivariant section is accordingly determined
by its values on the identity elements of $G$. Given any section $\zeta$ of
$\Bar{K}$, let $\zeta *\til0$ denote the section $(g,h)\mapsto
(\zeta(g),\til0_h)$ of $\Om*\Om\to G*G$. Then $\zeta$ is $G$-equivariant
iff $\zeta*\til0\proj{\kappa}\zeta$.

Now if $(\xi,X)$ is morphic, we have $\xi*\xi\proj{\kappa}\xi$, and so
$\Bar{\kappa}*\til0\proj{\kappa}\Bar{\kappa}$ implies that
$[\xi*\xi,\Bar{\kappa}*\til0]\proj{\kappa}[\xi,\Bar{\kappa}]$. Hence
$[\xi,\Bar{\kappa}]*\til0\proj{\kappa}[\xi,\Bar{\kappa}]$ and so
$[\xi,\Bar{\kappa}]$ is $G$-equivariant. Since $[\xi,\Bar{\kappa}]\circ 1 =
D_\xi(\kappa)$, we must have $[\xi,\Bar{\kappa}] = \Bar{D_\xi(\kappa)}$.

The last statement now follows immediately.
\boom

Thanks to Lemma \ref{lem:D} we can define
\begin{equation}                               \label{eq:bracket3}
[A(\xi),\kappa^\hs] = D_\xi(\kappa)^\hs
\end{equation}
for any star section $(\xi,X)$ and any $\kappa\in\Ga K$. By virtue of
Proposition~\ref{prop:generation}, we can now define a bracket on
$\Ga_{AG}A\Om$ by extending the conditions (\ref{eq:bracket1}),
(\ref{eq:bracket2}) and (\ref{eq:bracket3}) by the nonlinearity condition
\begin{equation}                               \label{eq:nonlin}
[\Xi,F\dtimes\Upsilon] = F\dtimes[\Xi,\Upsilon]\dpl{\bf a}(\Xi)(F)\dtimes
\Upsilon
\end{equation}
which must hold for all $\Xi,\Upsilon\in\Ga_{AG}(A\Om)$ and $F\in C(AG)$.

\begin{thm}                                  \label{thm:Om}
Let $(\Om;G,A;M)$ be an \LAgpd. Then the bundle $A\Om\to AG$ is a Lie
algebroid with respect to the anchor ${\bf a} = j_G^{-1}\circ A(\tila)\colon
A\Om\to TAG$ and the bracket just defined.
\end{thm}

We call this structure the {\em prolonged Lie algebroid structure} induced
from $\Om\to G$.

For the proof, we first need a Lemma, which involves the calculus of star
vector fields developed in \cite{MackenzieX:1998}.

\begin{lem}
For $(\xi,X)$ a star section of $\Om\to G$ with $\tila(\xi) = x$, and
$\kappa\in\Ga K$, we have
$$
{\bf a}\circ A(\xi) = \Tilde{x},\qquad
{\bf a}\circ\kappa^\hs = \da_{AG}(\kappa)^\upa,\qquad
\da_{AG}(D_\xi(\kappa)) = D_x(\da_{AG}(\kappa)).
$$
\end{lem}

\pf
Firstly, we have ${\bf a}\circ A(\xi)
= j_G^{-1}\circ A(\tila\circ\xi) = j_G^{-1}\circ A(x) = \Tilde{x},$ as
stated. The second identity is proved similarly.

For the third equation, we can calculate the LHS by finding
$\tila([\xi,\Bar{\kappa}](1_m))$, since we know that $D_\xi(\kappa)(m)$ is
a core element of $\Om$, and $\da_{AG}$ is the core of
$\tila\colon\Om\to TG.$ Now $\tila$ is the anchor of a Lie algebroid
structure, so
$$
\tila([\xi,\Bar{\kappa}](1_m)) = [x,\tila(\Bar{\kappa})](1_m).
$$
Write $W = \da_{AG}(\kappa)\in\Ga AG.$ Then $\tila(\Bar{\kappa}) =
\Vec{W}$, the right-invariant vector field on $G$ corresponding to $W$.
Finally, $[x,\Vec{W}](1_m) = D_x(W)(m)$.
\boom

\noindent{\sc Proof of Theorem \ref{thm:Om}.}\
It follows from the Lemma that
${\bf a}([\Xi,\Upsilon]) = [{\bf a}(\Xi),{\bf a}(\Upsilon)]$
for sections of the form $A(\xi)$ or $\kappa^\hs$. Likewise, the Jacobi
identity is easily checked for sections of these forms.

We must check that (\ref{eq:nonlin}) is consistent with (\ref{eq:bracket1}),
(\ref{eq:bracket2}), and (\ref{eq:bracket3}) in the cases where they overlap.
To verify consistency between (\ref{eq:bracket3}) and (\ref{eq:nonlin}) for
$\Xi = A(\xi),\ F = f\circ q_{AG},\ \Upsilon = \kappa^\hs$, we have
\begin{eqnarray*}
[A(\xi),(f\circ q_{AG})\dtimes\kappa^\hs] & = & [A(\xi), (f\kappa)^\hs]\\
     & = & D_\xi(f\kappa)^\hs\\
     & = & (fD_\xi(\kappa) + a(X)(f)\kappa)^\hs\\
     & = & (f\circ q_{AG})\dtimes D_\xi(\kappa)^\hs \dpl
                       (a(X)(f)\circ q_{AG})\dtimes \kappa^\hs\\
     & = & (f\circ q_{AG})\dtimes[A(\xi),\kappa^\hs] \dpl
                       (a(X)(f)\circ q_{AG})\dtimes \kappa^\hs
\end{eqnarray*}
and it remains to prove that $a(X)(f)\circ q_{AG} =
{\bf a}(A(\xi))(f\circ q_{AG})$. But ${\bf a}(A(\xi)) = \Tilde{W}$, where
$W = \Tilde{a}(\xi)$ is a star vector field on $G$ lying above $a(X)$, and
so the result follows.

Now taking $\Xi = \lambda^\hs,\ F = f\circ q_{AG},\ \Upsilon = \kappa^\hs$,
it is easy to see that both sides of (\ref{eq:nonlin}) are zero.

The only other relation requiring a check of consistency is
$$
A((f\circ\alpha)\xi) = (f\circ q_{AG})\dtimes A(\xi),
$$
where $\xi$ is a star section, and $f\in C(M)$, and the proof is very similar
to the preceding cases.

Thus the bracket is consistently defined, and the extension to the general
case of ${\bf a}([\Xi,\Upsilon]) = [{\bf a}(\Xi),{\bf a}(\Upsilon)]$ and the
Jacobi identity now follow from the nonlinearity condition (\ref{eq:nonlin}).
\boom

\begin{df}                            \label{df:dlalagpd}
The double vector bundle (\ref{eq:AOm}) together with the given Lie algebroid
structures on $A\to M,\ AG\to M$ and $A\Om\to A$ and the Lie algebroid
structure on $A\Om\to AG$ just defined, is the {\em double Lie algebroid of
the \LAgpd\ $(\Om:G,A,M)$.}
\end{df}

To justify this definition, we will in a future paper define an abstract
concept of double Lie algebroid and verify that (\ref{eq:AOm}) satisfies it.
In Example \ref{ex:plg} below we show how Theorem~\ref{thm:Om} includes the
construction of the Lie bialgebra of a Poisson Lie group (for example
\cite{LuW:1990}) and the construction of the Lie bialgebroid of a Poisson
groupoid \cite{MackenzieX:1994}. Further examples of different types are given
in \S\ref{sect:dladlg}. We first need the most basic example.

\begin{ex}\rm               \label{ex:first}
Let $A$ be any Lie algebroid on base $M$, and give $\Om = A\times A$ the
\LAgpd\ structure $(A\times A; M\times M, A; M)$ described at the end of
\cite[4.4]{Mackenzie:1992}. Then the Lie algebroid structure on $A\Om = TA$
with base $AG = TM$ defined above coincides with that defined on $TA\to TM$
in \cite[5.1]{MackenzieX:1994}.
\end{ex}

\begin{ex}\rm                \label{ex:plg}
Any Poisson Lie group $G$ gives rise to an \LAgpd\
$(T^*G;G,\frak{g}^*;\{1\})$, as described in \cite[4.12]{Mackenzie:1992}.
More generally, any Poisson groupoid $G\gpd P$ gives rise to an \LAgpd\
$(T^*G;G,A^*G;P)$ with the structures described in
\cite[\S8]{MackenzieX:1994}. We now show that the prolonged Lie algebroid
structure on $AT^*G\to AG$ admits a second natural description in terms of
the calculus of \cite{MackenzieX:1998}.

Consider at first a Lie groupoid $G\gpd P$. Then $T^*G\gpd A^*G$ has a
natural structure of symplectic groupoid and there is consequently an
isomorphism  $s\colon T^*A^*G\to AT^*G$ of Lie algebroids over $A^*G$,
where $T^*A^*G$ has the cotangent Lie algebroid structure for the Poisson
structure on $A^*G$, and $AT^*G$ is the Lie algebroid of $T^*G\gpd A^*G$
\cite{CDW}. It was shown in \cite[7.3]{MackenzieX:1994} that this isomorphism has an
alternative description. Firstly, by taking the dual of $j_G\colon TAG\to
ATG$ over $AG$, and identifying the dual of $ATG\to AG$ with $AT^*G\to AG$
via the induced pairing, we obtain a map $j'_G\colon AT^*G\to T^*AG$ which
is an isomorphism of double vector bundles over $AG$ and $A^*G$. Composing
this with the canonical isomorphism $R\colon T^*A^*G\to T^*AG$ from
\cite[5.5]{MackenzieX:1994}, we have two isomorphisms of double vector bundles, shown
in Figure~\ref{fig:Rs}.
\begin{figure}[htb]
\begin{picture}(350,280)(-100,0)
\put(-70,205){$\matrix{&&q_{T^*G}&\cr
                      &AT^*G&\lrah &A^*G\cr
                      &&&\cr
                A(c_G)&\Bigg\downarrow& &\cr
                      &&&\cr
                      &AG& & \cr}$}
\put(-10,220){\vector(1,-1){60}}            \put(30,200){$s^{-1}$}
\put(20,120){$\matrix{&&c_{A^*G}&\cr
                     &T^*A^*G&\lrah &A^*G\cr
                     &&&\cr
                     &\Bigg\downarrow&&\cr
                     &&&\cr
                     &AG& &  \cr}$}
\put(80,130){\vector(1,-1){60}}          \put(130,90){$R$}
\put(100,30){$\matrix{&&   &\cr
                     &T^*AG&\lrah &A^*G\cr
                     &&&\cr
               c_{AG}&\Bigg\downarrow&&\cr
                     &&&\cr
                     &AG& & \cr}$}
\end{picture}
\caption{\ \label{fig:Rs}}
\end{figure}
By \cite[7.3]{MackenzieX:1994}, $R\circ s^{-1} = j'_G$. In Figure~\ref{fig:Rs}
we have omitted the base $P$ and identity arrows for clarity.

Now suppose that $G\gpd P$ is a Poisson groupoid. Then $AT^*G\to AG$ has the
prolonged Lie algebroid structure (\ref{thm:Om}) from the cotangent structure
on $T^*G\to G$. On the other hand, $T^*AG\to AG$ has the cotangent Lie
algebroid structure from the Poisson structure on $AG$ dual to the Lie
algebroid structure on $A^*G$ \cite{Weinstein:1988}.

\begin{thm}[{\cite[\S7]{MackenzieX:1998}}]
For any Poisson groupoid $G\gpd P$, $j'_G\colon AT^*G\to T^*AG$ is an
isomorphism of Lie algebroids over $AG$ with respect to the structures
just described.
\end{thm}

Figure~\ref{fig:Rs} can now be regarded as depicting isomorphisms between
three double Lie algebroids. Firstly, $(AT^*G;AG,A^*G;P)$ is the double
Lie algebroid of $(T^*G;G,A^*G;P)$ and has the structures of
\ref{df:dlalagpd}. Secondly, $T^*A^*G\to A^*G$ is the cotangent Lie
algebroid of the Poisson structure on $A^*G$ dual to the Lie algebroid
structure of $AG$. The structure on $T^*A^*G\to AG$ is either transported
from the cotangent structure of $T^*AG\to AG$ via $R$, or is transported
from the prolonged structure on $AT^*G\to AG$ via $s$. The structures on
$T^*AG$ may be described similarly or obtained from the self-dual nature
of Lie bialgebroids \cite{MackenzieX:1994}.

In the case of Poisson Lie groups these structures may also be defined in
terms of the coadjoint actions and dressing transformations.
\end{ex}

\section{\bf THE DOUBLE LIE ALGEBROID OF A DOUBLE LIE GROUPOID}
\label{sect:dladlg}

Throughout this section we consider a double Lie groupoid $(S;H,V;M)$,
as defined in \cite[Definition 2.1]{Mackenzie:1992}. Applying the Lie
functor to the vertical groupoid structure, as in
\cite[\S4]{Mackenzie:1992}, we obtain the vertical \LAgpd\
$(A_VS;H,AV;M)$;
\begin{figure}[htb]
\begin{picture}(395,110)
\put(130,50){$\matrix{&&A(\tilalpha_H),A(\tilbeta_H)&&\cr
        &A_VS  &\lgpd&AV&\cr
        &&&&\cr
\tilq_V &\Bigg\downarrow&&\Bigg\downarrow&\cr
        &&&&\cr
        &H  &\lgpd&M,&\cr
        &&\alpha_H,\beta_H&&\cr}$}
\put(200,50){$TH\stackrel{T(\alpha_H),T(\beta_H)}{\lgpd}\ TM$}
\put(175,75){\vector(1,-1){15}}    \put(190,65){$\tila_V$}
\put(290,75){\vector(1,-1){15}}           
\put(195,45){\vector(-1,-1){15}}
\put(305,45){\vector(-1,-1){15}}
\end{picture}         \caption{\ \label{fig:AVS}}
\end{figure}
see Figure~\ref{fig:AVS}. The construction of \S\ref{sect:dlalag} now
leads to a double vector bundle
\begin{equation}                               \label{eq:A2S}
\matrix{&&\ts{q}_H&&\cr
        &A^2S&\lrah &AV&\cr
        &&&&\cr
A(\tilq_V)&\Bigg\downarrow& &\Bigg\downarrow&q_V\cr
        &&&&\cr
        &AH&\lrah &M&\cr
        &&q_H&&\cr}
\end{equation}
where $A^2S = A(A_VS)$ has two Lie algebroid structures: over base $AV$ the
structure of the Lie algebroid of the Lie groupoid $A_VS\gpd AV$, and over
$AH$ the prolongation of the Lie algebroid structure on $A_VS\to H$.

Equally well, we can first take the horizontal \LAgpd\ $(A_HS;AH,V;M)$ as
in Figure~\ref{fig:hLAgpd}.
\begin{figure}[hbt]
\begin{picture}(395,110)
\put(70,50){$\matrix{      &&\tilq_H &&\cr
                            &A_HS    &\vlra &V &\cr
                            &&&&\cr
A(\tilalpha_V),A(\tilbeta_V)&\vgpd&&\vgpd&\cr
                            &&&&\cr
                            &AH &\vlra &M,&\cr
                            &&&&\cr}$}
\end{picture}         \caption{\ \label{fig:hLAgpd}}
\end{figure}

\noindent
Applying the construction of \S\ref{sect:dlalag} to $A_HS$ yields a double
vector bundle $A_2S = A(A_HS)$
\begin{equation}                               \label{eq:A_2S}
\matrix{&&A(\tilq_H)&&\cr
        &A_2S&\lrah &AV&\cr
        &&&&\cr
\ts{q}_V&\Bigg\downarrow& &\Bigg\downarrow&q_V\cr
        &&&&\cr
        &AH&\lrah &M&\cr
        &&q_H&&\cr}
\end{equation}
and two Lie algebroid structures on $A_2S$: that over base $AH$ being the Lie
algebroid of $A_HS\gpd AH$, and that over $AV$ being the prolongation of
$A_HS\to V$.

Before proceeding, we need to consider a basic example.

\begin{ex}\rm                \label{ex:pairs}
Let $G$ be a Lie groupoid on $M$, and let $S = G\times G$ have the double
Lie groupoid structure of \cite[2.3]{Mackenzie:1992} with $H = M\times M$
and $V = G$. Then, as in \cite[4.4]{Mackenzie:1992}, the vertical and
horizontal \LAgpd s are respectively
\begin{equation}                    \label{eq:AGAGTG}
\matrix{&AG\times AG&\lgpd&AG&\cr
        &&&&\cr
        &\Bigg\downarrow&&\Bigg\downarrow&\cr
        &&&&\cr
        &M\times M&\lgpd&M,&\cr
        }
\qquad\mbox{and}\qquad
\matrix{                    &TG&\vlra &G&\cr
                            &&&&\cr
                            &\vgpd&&\vgpd&\cr
                            &&&&\cr
                            &TM &\vlra &M.&\cr
                            &&&&\cr}
\end{equation}
We thus obtain $A^2S = TAG$ and $A_2S = ATG$.

In \cite[7.1]{MackenzieX:1994}, it was shown that $j_G\colon TAG\to ATG$
is an isomorphism of double vector bundles from
$$
\matrix{&&p_{AG}&&\cr
        &TAG&\lrah&AG&\cr
        &&&&\cr
  T(q_G)&\Bigg\downarrow&&\Bigg\downarrow&\cr
        &&&&\cr
        &TM&\lrah&M&\cr
        &&&&\cr}
\mbox{\qquad to\qquad}
\matrix{&&A(p_G)&&\cr
        &ATG&\lrah&AG&\cr
        &&&&\cr
  q_{TG}&\Bigg\downarrow&&\Bigg\downarrow&\cr
        &&&&\cr
        &TM&\lrah&M&\cr
        &&&&\cr}
$$
preserving $AG$ and $TM$; $j_G$ is a restriction of the canonical involution
$J_G\colon T^2G\to T^2G$ associated with the manifold $G$. Further,
\cite[7.1]{MackenzieX:1994} shows that $j_G$ is an isomorphism of the Lie algebroid
structures over $TM$: it sends the tangent Lie algebroid structure of
$T(AG)\to TM$ (which by Example \ref{ex:first} is the prolongation of
$AG\times AG\to M\times M$) to the Lie algebroid structure $A(TG)\to TM$
arising from the tangent groupoid $TG\gpd TM$.

An element $\Xi\in TAG$ can be written as $\frac{\da^2 g}{\da t\da u}(0,0)$
where $g\colon\R^2\to G$ is a smooth map such that if
$m = \beta\circ g\colon\R^2\to M$ then $\alpha(g(t,u)) = m(t,0)$ and
$g(t,0) = 1_{m(t,0)}$ for all $t,u$. The notation is meant to indicate that
we first differentiate with respect to $u$ and obtain a curve in $AG$, whose
derivative with respect to $t$ is $\Xi$. Accordingly,
$p_{AG}(\Xi) = \frac{\da g}{\da u}(0,0), \
T(q_G)(\Xi) = \frac{\da m}{\da t}(0,0)$ and
$$
j_G\left(\frac{\da^2 g}{\da t\da u}(0,0)\right) =
\frac{\da^2 g}{\da u\da t}(0,0).
$$

It is easily seen that $j_G$ is also an isomorphism of the Lie algebroid
structures over $AG$. The structure in $ATG\to AG$ is defined by
$$
[A(\xi),A(\eta)] = A([\xi,\eta]),\qquad
[A(\xi),X^\hs] = D_\xi(X)^\hs,\qquad
[X^\hs,Y^\hs] = 0,
$$
where $\xi$ and $\eta$ are star vector fields on $G$ and $X,Y\in\Ga AG$.
Now from \cite[\S3]{MackenzieX:1998} we have
$$
j_G^{-1}\circ A(\xi) = \Tilde{\xi},\qquad
j_G^{-1}\circ X^\hs = X^\upa,
$$
and
\begin{equation}                      \label{eq:MX3}
[\Tilde{\xi},\Tilde{\eta}] = \Tilde{[\xi,\eta]},\qquad
[\Tilde{\xi},X^\upa] = D_\xi(X)^\upa,\qquad
[X^\upa,Y^\upa] = 0,
\end{equation}
from which the result follows.

Thus $j_G$ preserves both the horizontal and the vertical Lie algebroid
structures.

\end{ex}
\boom

We now show that a similar phenomenon holds for any double Lie groupoid
$(S;H,V;M)$.

\begin{thm}                                 \label{thm:jS1}
The canonical involution $J_S\colon T^2S\to T^2S$ of the manifold $S$
restricts to an isomorphism of double vector bundles $j_S\colon A^2S\to A_2S$
which preserves $AH$ and $AV$.
\end{thm}

\pf
We regard elements of $T^2S$ as second derivatives
$$
\Xi = \frac{\da^2s}{\da t\da u}(0,0)
$$
where $s\colon\R^2\to S$ is a smooth square of elements, and the notation
means that $s$ is first differentiated with respect to $u$, yielding a
curve $\xi_t = \frac{\da s}{\da u}(t,0)$ in $TS$ with
$\ddt{\xi_t} = \Xi$. Thus
$$
\frac{\da s}{\da u}(0,0) = p_{TS}(\Xi),\qquad
\frac{\da s}{\da t}(0,0) = T(p_S)(\Xi),
$$
and $J_S(\Xi) = \frac{\da^2s}{\da u\da t}(0,0).$

If $\Xi\in A^2S$ then $s$ can be chosen to be of the form
$$
\matrix{                  &\raise.7ex\hbox{$h(t,u)$}&\cr
\raise.35in\hbox{$v(t,u)$}&\ssq&\raise.35in\hbox{$v(0,u)$}\cr
                          &h(t,0)&\cr}\hskip-1.2in\raise5pt\hbox{$s(t,u)$}
$$
where $h\colon\R^2\to H,\ v\colon\R^2\to V$ are smooth and satisfy the source
and target conditions implicit in the diagram, and where
$s(t,0) = \tilone^V_{h(t,0)},\ s(0,u) = \tilone^H_{v(0,u)}.$
We now have $\xi_t = \frac{\da s}{\da u}(t,0)\in A_VS$ with
$$
\tilq_V(\xi_t) = h(t,0),\quad
A(\tilbeta_V)(\xi_t) = \frac{\da v}{\da u}(t,0),\quad
A(\tilalpha_V)(\xi_t) = \frac{\da v}{\da u}(0,0),
$$
and $\xi_0 = A(\tilone^H)(\frac{\da v}{\da u}(0,0)).$ Consequently
$$
A(\tilq_V)(\Xi) = \frac{\da h}{\da t}(0,0),\qquad
\ts{q}_H(\Xi) = \frac{\da v}{\da u}(0,0).
$$

Interchanging the roles of $t$ and $u$, it is easily seen that
$j_S(\Xi)\in A_2S$. Since $j_S$ is a restriction of $J_S$ it is injective
and, by a dimension count, an isomorphism of double vector bundles.
\boom

Theorem \ref{thm:jS1} can also be proved by the diagram chasing methods
used in \cite[Proposition 1.5]{Mackenzie:1992} and
\cite[7.1]{MackenzieX:1994}.

\begin{thm}                             \label{thm:jS2}
Let $(S;H,V;M)$ be a double Lie groupoid. Then $j_S$, regarded as a morphism
of vector bundles over $AH$, is an isomorphism of Lie algebroids from the
prolonged structure on $A(A_VS)\to AH$ to the Lie algebroid of $A_HS\gpd AH$.
\end{thm}

\pf
That $j_S$ commutes with the anchors is proved by the same methods as in
Theorem \ref{thm:jS1}. It must be shown that for $\Xi\in A^2S$, we have
$j_H\circ\ts{a}_V\circ j_S(\Xi) = A(\tila_V)(\Xi)$. Taking $\Xi =
\displaystyle\frac{\da^2s}{\da t\da u}(0,0)$ as above, both sides are equal
to $\displaystyle\frac{\da^2h}{\da t\da u}(0,0)$.

The remainder of the proof depends on a series of subsidiary results,
\ref{lem:tilRi} to \ref{cor:D}.

Given a star section $(\xi,X)$ of $A_VS$, write $\Tilde{\xi} =
j_S\circ A(\xi)$. As with sections of the Lie algebroid of any Lie groupoid,
$\xi$ induces a vector field $\Vec{\xi}$ on $S$, which is right-invariant
with respect to the groupoid structure on base $H$. Similarly, $X$ induces
a right invariant vector field $\Vec{X}$ on $V$. Since $(\xi,X)$ is a star
section, it follows that $(\Vec{\xi},\Vec{X})$ is a star vector field on the
groupoid $S\gpd V$. Hence, by the constructions of \cite{MackenzieX:1998}, it induces
a linear vector field $((\Vec{\xi})\twiddle,\Vec{X})$ on $A_HS$, here regarded
as a Lie algebroid over $V$.

On the other hand, $\Tilde{\xi}$ is a section of the Lie algebroid of
$A_HS\gpd AH$, and therefore induces a right invariant vector field
$(\Tilde{\xi})\riri$ on $A_HS$.

\begin{lem}                               \label{lem:tilRi}
For any star section $(\xi,X)$ of the \LAgpd\ $A_VS$,
$(\Tilde{\xi})\riri = (\Vec{\xi})\twiddle.$
\end{lem}

\pf
Let $\Exp t\xi$ be the exponential for $\xi\in\Ga(A_VS)$. We follow the
conventions of \cite{Mackenzie:LGLADG}, and assume for convenience that
$\Exp t\xi$ is global. Similarly, let $\Exp tX$ be the exponential for
$X\in\Ga AV$. Since $\tilalpha_H\circ\xi = X\circ\alpha_H$ and
$\xi\circ 1^H = \tilone^H\circ X$, it follows that
$$
\tilalpha_H\circ\Exp t\xi = \Exp tX\circ\alpha_H,\qquad
(\Exp t\xi)\circ 1^H = \tilone^H\circ\Exp tX.
$$

The corresponding flow for $\Vec{\xi}$ is $L_{\Exp t\xi}\colon S\to S,\
s\mapsto \Exp t\xi(\tilbeta_V s)\,\vcomp\, s$ and is a star map over
$L_{\Exp tX}\colon V\to V$; that is;
\begin{equation}                            \label{eq:flows}
\tilalpha_H\circ L_{\Exp t\xi} = L_{\Exp tX}\circ\alpha_H
\qquad\mbox{and}\qquad
L_{\Exp t\xi}\circ \tilone^H = \tilone^H\circ L_{\Exp tX}.
\end{equation}
We can therefore apply the Lie functor and get $A(L_{\Exp t\xi})\colon
A_HS\to A_HS$, which is a linear map over $L_{\Exp tX}\colon V\to V$.
By \cite[\S3]{MackenzieX:1998}, this is the flow of $(\Vec{\xi})\twiddle$.

Now it is straightforward to check that $A(L_{\Exp t\xi})$ is also the flow
of $(\Tilde{\xi})\riri$, and from this the result follows.
\boom

\begin{cor}                                     \label{cor:Tildes}
For any star sections $(\xi,X)$ and $(\eta,Y)$ of the \LAgpd\ $A_VS$,
$[\Tilde{\xi},\Tilde{\eta}] = \Tilde{[\xi,\eta]}.$
\end{cor}

\pf
As for any Lie groupoid, $[\Tilde{\xi},\Tilde{\eta}]\riri =
[(\Tilde{\xi})\riri,(\Tilde{\eta})\riri] =
[(\Vec{\xi})\twiddle,(\Vec{\eta})\twiddle]$, by \ref{lem:tilRi}. Now this is
$[\Vec{\xi},\Vec{\eta}]\twiddle$ by (\ref{eq:MX3}),
and using $[\Vec{\xi},\Vec{\eta}] = [\xi,\eta]\riri$ and \ref{lem:tilRi}
again, the result is clear.
\boom

Next consider $\kappa\in\Ga K$, where $K$ is the core of $A_VS$. From
(\ref{eq:coresection}) and (\ref{eq:coreds}), we have $\kappa^\hs$, a section
of $A^2S\to AH$ and $\kappa^\ds$, a section of $A^2S\to AV$. On the other
hand, we can also consider $K$ as the core of $A_HS$, and we denote the
corresponding sections respectively by $\kappa^\cls$, a section of
$A_2S\to AV$, and $\kappa^\spades$, a section of $A_2S\to AH$. This device
depends on the fact that the cores of $A_HS$ and $A_VS$ can be canonically
identified (Proposition \ref{prop:cores}). It is now easily checked that
\begin{equation}               \label{eq:hsjss}
j_S\circ\kappa^\hs = \kappa^\spades, \qquad
j_S\circ\kappa^\ds = \kappa^\cls.
\end{equation}

In the next result, $\BarH{\kappa}$ is the section of $A_HS\to V$ induced by
$\kappa$. As a section of the vector bundle $A_HS$, it induces a vertical
vector field on $A_HS$. On the other hand, $\kappa^\spades$ is a section of
the Lie algebroid of the groupoid $A_HS\gpd AH$, and consequently induces
a right invariant vector field on $A_HS$.

\begin{lem}                    \label{lem:upri}
For any $\kappa\in\Ga K$,
$(\BarH{\kappa})^\upa = (\kappa^\spades)\riri.$
\end{lem}

\pf
Take any $\zeta\in A_HS\vert_v$ with $A(\tilbeta_V)(\zeta) = y\in A_nH$. We
know from (\ref{eq:coreds}) that
$$
\kappa^\spades(y) = \ddt{(A(\tilone^V)(y) + t\kappa(n))}
$$
where $+$ is the addition in $A_HS\to V$. Right translating this curve by
$\zeta$, we get
$$
(A(\tilone^V)(y) + t\kappa(n))\zeta =
(A(\tilone^V)(y) + t\kappa(n))(\zeta + \til0_v) =
A(\tilone^V)(y)\zeta + t(\kappa(n)\til0_v)
$$
by the interchange laws. This in turn is equal to
$\zeta + t\BarH{\kappa}(v)$, and
$$
(\BarH{\kappa})^\upa(\zeta) = \ddt{(\zeta + t\BarH{\kappa}(v))}
$$
is just the definition of the vertical lift.
\boom

In the next result, the $D$ operator on the RHS refers to the \LAgpd\ $A_VS$
and is as defined in (\ref{eq:D}), whereas the $D$ on the LHS refers to the
groupoid $S\gpd V$, and is defined in \cite[\S3]{MackenzieX:1998}.

\begin{prop}
Let $(\xi,X)$ be a star section of $A_VS$, and let $\kappa\in\Ga K$. Then
$$
D_{\Vec{\xi}}(\BarH{\kappa}) = \BarH{D_\xi(\kappa)}.
$$
\end{prop}

\pf
We must prove that for all $v\in V$,
$$
[\Vec{\xi}, \VecH{\kappa}](\tilone^H_v) =
   [\xi,\BarV{\kappa}](1^H_{\beta v})\til0_v.
$$
On the LHS, $\Vec{\xi}$ is the right invariant vector field on $S\gpd H$
corresponding to $\xi$, here considered as a star vector field on $S\gpd V$,
whilst $\VecH{\kappa}$ is the right invariant vector field, with respect to
the horizontal groupoid structure on $S$, corresponding to
$\BarH{\kappa}\in\Ga_V(A_HS)$. On the RHS we have the bracket of two sections
of the Lie algebroid $A_VS$, and the multiplication in $A_HS\gpd AH$.

Assume as in the proof of Lemma \ref{lem:tilRi} that $\Exp t\xi$ is a
global exponential for $\xi$. Write $m = \beta v$ and
$m_t = \beta\Exp tX(m)$. Then, using (\ref{eq:flows}),
$$
[\Vec{\xi}, \VecH{\kappa}](\tilone^H_v)
= -\ddt{(L_{\Exp t\xi})_*(\VecH{\kappa})(\tilone^H_v)}\\
= -\ddt{T(L_{\Exp t\xi})(\kappa(m_{-t})\til0_{\Exp -tX(m)v})}.
$$
On the other hand,
$[\xi,\BarV{\kappa}](1^H_{\beta v}) =
         -\ddt{T(I_{\Exp t\xi})(\kappa(m_{-t}))}$
where $I$ is the conjugation in $S\gpd H$;
$$
I_{\Exp t\xi}(s) =
\Exp t\xi(\tilbeta_V(s))\,\vcomp\, s\,\vcomp\,
\left(\Exp t\xi(\tilalpha_V(s))\right)^{-V}.
$$

If we now write
$\kappa(m_{-t}) = \left.\frac{\da}{\da s}c(s,t)\right|_{s = 0}$
where $c(s,t)$ is a curve in $C$ in the $\alpha$-fiber above $m_{-t}$, then
$$
T(I_{\Exp t\xi})(\kappa(m_{-t})) =
\left.\frac{\da}{\da s}\Exp t\xi(\tilbeta_V(c(s,t)))\,\vcomp\,c(s,t)\,\vcomp\,
\tilone^H_{\Exp tX(m_{-t})^{-1}}\right|_{s = 0}.
$$
Recalling that $(\Exp tX(m_{-t}))^{-1} = \Exp -tX(m)$, it now follows that
$$
T(L_{\Exp t\xi})(\kappa(m_{-t})\til0_{\Exp -tX(m)v})
= T(I_{\Exp t\xi})(\kappa(m_{-t}))\til0_v,
$$
and this completes the proof.
\boom

\begin{cor}                                \label{cor:D}
Let $(\xi,X)$ be a star section of $A_VS$, and let $\kappa\in\Ga K$. Then,
in the Lie algebroid of $A_HS\gpd AH$,
$$
[\Tilde{\xi},\kappa^\spades] = D_\xi(\kappa)^\spades.
$$
\end{cor}

\pf
It suffices to prove that the corresponding right invariant vector fields
on $A_HS$ are equal. Now, using Lemmas \ref{lem:tilRi} and \ref{lem:upri},
the right invariant vector field corresponding to the LHS is
$[(\Vec{\xi})\twiddle, (\BarH{\kappa})^\upa]$.
Here $(\Vec{\xi},\Vec{X})$ is a star vector field on $S\gpd V$, and so induces
a vector field on $A_HS$. Similarly, $(\BarH{\kappa})^\upa$ is the vertical
lift to $A_HS$ of a section of the Lie algebroid $A_HS$. Applying
(\ref{eq:MX3}), we therefore get $D_{\Vec{\xi}}(\BarH{\kappa})^\upa$. The
result now follows from using Lemma \ref{lem:upri} again.
\boom

We can now complete the proof of Theorem \ref{thm:jS2}. By Proposition
\ref{prop:generation} we know that $\Ga_{AH}(A^2S)$ is generated by sections
of the form $A(\xi)$, where $(\xi,X)$ is a star section of $A_VS$, and those
of the form $\kappa^\hs$, where $\kappa\in\Ga K$. Since we have
$j_S\circ A(\xi) = \Tilde{\xi}$ and $j_S\circ\kappa^\hs = \kappa^\spades$,
the equations from \ref{cor:Tildes}, \ref{cor:D} and \ref{prop:upa} prove
that $j_S$ preserves the brackets.
\boom

Theorem \ref{thm:jS2} clearly includes the isomorphisms of Example
\ref{ex:pairs}. The isomorphisms of Example \ref{ex:plg}, however, are
of a different type.

It is now clear that there are two, canonically isomorphic, double Lie
algebroids associated with a double Lie groupoid. We make an arbitrary
choice.

\begin{df}
The double vector bundle $A^2S$ of (\ref{eq:A2S}), equipped with the Lie
algebroid structure on $A^2S\to AV$ of the Lie groupoid $A_VS\gpd AV$ and the
prolonged Lie algebroid structure on $A^2S\to AH$ from $A_VS\to H$,
is the {\em double Lie algebroid of the double Lie groupoid $(S;H,V;M)$.}
\end{df}

\begin{ex}\rm
Take $G = M\times M$ in Example~\ref{ex:pairs}, where $M$ is any manifold.
Thus $S = M^4$ consists of quadruples of points from $M$. From \ref{ex:pairs}
we obtain $A^2S = T^2M$ and $A_2S = T^2M$ with $j_S = J_M$ the canonical
involution of $M$. The prolonged structure on $T(p)\colon T^2M\to TM$ has
anchor $J = J_M$ and bracket $[\xi,\eta] = J[J\xi,J\eta]$ for all sections
$\xi,\eta$ of $T(p)$.
\end{ex}

\begin{ex}\rm                              \label{ex:quads}
Let $H$ and $V$ be Lie groupoids on the same base $M$, and let $S$ be the
double Lie groupoid $\tsq\,(H,V)$ of \cite[2.4]{Mackenzie:1992}; assume that
the anchors of $H$ and $V$ are suitably transversal.

Then, recalling \cite[4.5]{Mackenzie:1992}, the vertical \LAgpd\ $A_VS$ is
\begin{equation}
\raisebox{30pt}{$\chigh_H^{!!}(AV\times AV)\approx$\ }
\matrix{a_V^{!!}(TH)&\lgpd&AV&\cr
        &&&\cr
        \Bigg\downarrow&&\Bigg\downarrow&\cr
        &&&\cr
        H        &\lgpd&M.&\cr
        }
\end{equation}
Here $\chigh_H^{!!}(AV\times AV)$ is the pullback Lie algebroid of the
cartesian square $AV\times AV$ across the anchor
$\chigh_H\colon H\to M\times M$, and $a_V^{!!}(TH)$ is the pullback Lie
groupoid of $TH\gpd TM$ across $a_V\colon AV\to TM$.

Hence the horizontal structure on $A^2S = A(A_VS)$ is the pullback Lie
algebroid $a_V^{!!}(ATH)$ of $ATH\to TM$ across $a_V$. As a manifold
$$
A^2S = TAV\pback{T^2M}ATH,
$$
the pullback of $T(a_V)$ and $a_{TH}$. The anchor $\ts{a}_H\colon A^2S\to TAV$
is the natural map $TAV\pback{T^2M}ATH\to TAV$ and the bracket is of the
standard form for pullbacks: see \cite[\S1]{HigginsM:1990a}.

Performing the horizontal and vertical differentiations in reverse order, we
have
$$
A_2S = ATV\pback{T^2M}TAH,
$$
the pullback of $a_{TV}$ and $T(a_H)$. Over base $AH$ the Lie algebroid
structure on $A_2S$ is the pullback of $ATV\to TM$ over $a_H\colon AH\to TM$.
The canonical map $j_S\colon A^2S\to A_2S$ is $j_V\pback{T^2M}j_H^{-1}$.
Theorem \ref{thm:jS2} now shows that the prolonged structure on
$TAV\pback{T^2M}ATH\to AH$ is isomorphic to the pullback of $ATV\to TM$
across $a_H$.
\end{ex}

\begin{ex}\rm                                             \label{ex:comma}
Let $H$ and $V$ be Lie groupoids on base $M$, and let $\phi \colon H\to V$
be a base-preserving morphism. Let $\Theta = \Theta(H,\phi,V)$ be the comma
double groupoid of \cite[2.5]{Mackenzie:1992}. Then, recalling
\cite[4.7]{Mackenzie:1992}, the vertical \LAgpd\ $A_V\Th$ is
\begin{equation}
\raisebox{38pt}{$\alpha_H^{!!}(AV) \approx$}
\matrix{TH\sdp a_V&\lgpd &AV&\cr
        &&&\cr
        \Bigg\downarrow&&\Bigg\downarrow&\cr
        &&&\cr
        H&\lgpd &M.&\cr
        &&&\cr}
\end{equation}
Here $\alpha_H^{!!}(AV)$ is the pullback Lie algebroid of $AV$ across
$\alpha_H\colon H\to M$, and $TH\sdp a_V$ is the action groupoid arising
from the action of $TH\gpd TM$ on $a_V$ described in
\cite[4.7]{Mackenzie:1992}.

To describe the induced infinitesimal action of $ATH$ on $a_V\colon AV\to TM$,
recall from \cite{MackenzieX:1998} that in any Lie algebroid $A\to M$, a
section $X$ of
$A$ induces the vertical lift $X^\upa$ on $A$ and also a linear vector field
$(\Tilde{X},a(X))$ characterized by $\Tilde{X}(\ell_\psi) = \ell_{L_X(\psi)}$
for all sections $\psi$ of $A^*$.

The sections of $ATH\to TM$ are generated by those of the form $j_G\circ T(X)$
and those of the form $j_G\circ\Hat{X}$, where $X\in\Ga AH$ (see
\cite[7.1]{MackenzieX:1994}). Now it is easy to check that the infinitesimal
action is characterized by the equations
\begin{equation}
(j_G\circ T(X))^\dagger = \Tilde{\phi(X)},\qquad
(j_G\circ \Hat{X})^\dagger = \phi(X)^\upa.
\end{equation}

Since the Lie algebroid of an action Lie groupoid is the corresponding action
Lie algebroid, it follows that $A^2\Th\to AV$ is the action Lie algebroid
$ATH\sdp a_V$. As a manifold, $A^2\Th$ is $ATH\pback{TM}AV$, the pullback
of $q_{TH}$ and $a_V$.

On the other hand the horizontal \LAgpd\ of $\Th$ is
\begin{equation}
\raisebox{38pt}{$q_H^{!!}(V)\approx$}
\matrix{(AH\times AH)\sdp \chigh_V&\vlra &V&\cr
        &&&\cr
        \vgpd&&\vgpd&\cr
        &&&\cr
        AH&\vlra &M.&\cr
        &&&\cr}
\end{equation}
Here $q_H^{!!}(V)$ is the pullback groupoid of $V$ over $q_H\colon AH\to M$
and $(AH\times AH)\sdp\chigh_V$ is the action Lie algebroid arising from
the action of $AH\times AH$ on $\chigh_V$ given in
\cite[4.7]{Mackenzie:1992}.

Accordingly, $A_2\Th\to AH$ is the pullback Lie algebroid $q_H^{!!}(AV)$.
As a manifold, $A_2\Th$ is $TAH\pback{TM}AV$, and $j_{\Th}$ is
$j_H^{-1}\pback{TM}\id$.
\end{ex}
\boom

Notice throughout Examples \ref{ex:pairs}, \ref{ex:quads} and \ref{ex:comma},
that if the horizontal (say) structure is a pair groupoid, or a pullback, or
an action groupoid, then the corresponding property is inherited by the
double Lie algebroid. This follows from a version of
\cite[Lemma 4.6]{Mackenzie:1992} formulated for \LAgpd s and their double
Lie algebroids.

\begin{ex}\rm (Compare \cite[\S4.5]{Weinstein:1988}.)\label{ex:symp}
Let $\Ga\gpd P$ be any symplectic groupoid, and consider the double groupoid
$S = \Ga\times\Ga$ of Example \ref{ex:pairs}. Then the \LAgpd s
(\ref{eq:AGAGTG}) can be given as
$$
\matrix{&&      &&\cr
        &T^*P\times T^*P&\lgpd&T^*P&\cr
        &&&&\cr
        &\Bigg\downarrow&&\Bigg\downarrow&\cr
        &&&&\cr
        &P\times P&\lgpd&P&\cr
        &&&&\cr}
\mbox{\qquad and \qquad}
\matrix{&&      &&\cr
        &T\Ga&\lrah&\Ga&\cr
        &&&&\cr
        &\vgpd&&\vgpd&\cr
        &&&&\cr
        &TP&\lrah&P&\cr
        &&&&\cr}
$$
where the canonical isomorphism of Lie algebroids $T^*P\to A\Ga$ of \cite{CDW}
induces an isomorphism of \LAgpd s. Now $T\Ga\gpd TP$ is itself a symplectic
groupoid with respect to the tangent symplectic structure on $T\Ga$ and the
tangent Poisson structure on $TP$ and so $AT\Ga\isom T^*(TP)$. Thus the double
Lie algebroids are
\begin{equation}                       \label{eq:sympcase}
\matrix{&&      &&\cr
        &T(T^*P)&\lrah&T^*P&\cr
        &&&&\cr
        &\Bigg\downarrow&&\Bigg\downarrow&\cr
        &&&&\cr
        &TP&\lrah&P&\cr
        &&&&\cr}
\mbox{\qquad and\qquad}
\matrix{&&      &&\cr
        &T^*(TP)&\lrah&T^*P&\cr
        &&&&\cr
        &\Bigg\downarrow&&\Bigg\downarrow&\cr
        &&&&\cr
        &TP&\lrah&P&\cr
        &&&&\cr}
\end{equation}
with $j_S$ now the canonical isomorphism $T(T^*P)\to T^*(TP)$. Theorem
\ref{thm:jS2} applied to the vertical structures now shows that the
prolongation of the cotangent Lie algebroid structure on $T^*P$ is isomorphic
to the cotangent Lie algebroid of the tangent Poisson structure on $TP$.

Both double vector bundles in (\ref{eq:sympcase}) are duals of the double
tangent bundle $T^2P$ in the sense of Pradines' notion of dual
(see \cite[\S5]{MackenzieX:1994}). In these terms, the relations between the
Lie algebroid structures in (\ref{eq:sympcase}) reflect the fact that the
canonical involution $J_P\colon T^2P\to T^2P$ is a Poisson automorphism
with respect to the tangent of the tangent of the Poisson structure
on $P$ (see \cite[5.2]{Courant:1990t}).

The double Lie algebroid $T(T^*P)\isom T^*(TP)$ and these isomorphisms can
be defined for any (not necessarily integrable) Poisson manifold.
\end{ex}
\boom

Several further examples requiring more extensive developments, such as
the infinitesimal form of the theory associated with split double Lie
groupoids and crossed modules \cite{BrownM:1992} \cite[\S2]{Mackenzie:1992},
will be treated in other papers. The case of the Lie algebroid structures
associated with a vacant double Lie groupoid \cite[\S2]{Mackenzie:1992}
is dealt with by Mokri \cite{Mokri:1997}. In the final two sections of the
present paper we deal with the case of symplectic and Poisson double
groupoids. In the next section we consider the infinitesimal theory
associated with affinoids. \S\ref{sect:datcdg} and \S\ref{sect:sdgapdg}
are independent of \S\ref{sect:aff}.

\section{INFINITESIMAL STRUCTURES ASSOCIATED WITH AFFINOIDS}
\label{sect:aff}

The notion of affinoid structure takes several different forms: as a
ternary relation, and in the context of dual pairs in symplectic geometry,
it was introduced by Weinstein \cite{Weinstein:1990}; as a generalization
of principal bundles it was introduced by Kock \cite{Kock:1988} under the
name of {\em pregroupoid}; as a form of Morita equivalence for groupoids
it was introduced by Pradines \cite{Pradines:1977} as a {\em butterfly
diagram}. Forms of the notion, however, go back to the early part of the
century; see the references in \cite{Weinstein:1990} and \cite{Kock:1988}.

In \cite[\S3]{Mackenzie:1992} we gave proofs of the equivalences between
affinoid structures, butterfly diagrams and generalized principal bundles,
using simple functorial constructions from groupoid and double groupoid
theory. The key was to regard an affinoid structure as a type of double
groupoid, called {\em principal} in \cite[3.2]{Mackenzie:1992}. The
interest of affinoid structures for us here
is that the other two equivalent formulations, which are not overtly double
structures, provide a means for testing the correctness of the double Lie
algebroid construction. In this section we accordingly calculate the
infinitesimal invariants associated with affinoid structures, butterfly
diagrams and generalized principal bundles. Although the term \kwo{affinoid
structure} should strictly refer to the ternary relation, in this section we
will use it to mean the corresponding principal double Lie groupoid.

We first show that the double Lie algebroid of an affinoid structure
consists of a pair of conjugate flat partial connections adapted to
the two foliations of the affinoid (compare \cite[3.2]{Weinstein:1990}).
Providing the leaves of both foliations are simply--connected, such a
pair of partial connections is equivalent to the infinitesimal version
of butterfly diagram (\ref{thm:iasibd}, \ref{thm:ibdias}). However
extending the Atiyah sequence construction to generalized principal
bundles loses information and appears to be only of limited interest
(\ref{prop:gAs}).

Some of the results of this section (\ref{prop:diag}, \ref{lem:sum})
are special cases of results of Mokri \cite{Mokri:1997}; this reflects the
fact that an affinoid structure is a vacant double Lie groupoid. Here we can
give more direct proofs which are also of interest, however.

Consider an affinoid structure $S$ on a manifold $M$ with surjective
submersions $c\colon M\to Q_H$ and $b\colon M\to Q_V$. Orient the side
groupoids as $H = R(b)$ and $V = R(c)$ with quotient groupoids respectively
$G_h\gpd Q_H$ and $G_v\gpd Q_V$. In the notation of
\cite[\S3]{Mackenzie:1992}, $G_h = \tau^V_1(S)$ and $G_v = \tau^H_1(S)$.
Denote the fibrations defining the quotients by $\Tilde{c}\colon R(b)\to G_h$
(over $c$) and $\Tilde{b}\colon R(c)\to G_v$ (over $b$).

The morphism $\Tilde{b}\colon R(c)\to G_v$ defines a vertical subbundle
$T^{\Tilde{b}}(R(c)) \subseteq T(R(c)) = R(T(c)) = TM\pback{TQ_H}TM.$
On the other hand, applying the Lie functor to the morphism
$\Tilde{c}\colon R(b)\to G_h$ yields a Lie algebroid morphism
$A(\Tilde{c})\colon T^bM\to AG_h$.

\begin{lem}
$T^{\Tilde{b}}(R(c)) = R(A(\Tilde{c})).$
\end{lem}

\pf
Take $(Z,X)\in T^{\Tilde{b}}(R(c))$. Thus $Z\in T_z(M),\ X\in T_x(M)$ and
$$
T(c)(Z) = T(c)(X),\qquad\qquad T(\Tilde{b})(Z,X) = 0.
$$
Now $\Tilde{b}$ is a groupoid morphism over $b\colon M\to Q_V$ and so, taking
tangents,
$$
\dvb{T(R(c))}{TM}{TG_v}{TQ_V}{T(\Tilde{b})}{T(\alpha_V)}{T(\alpha)}{T(b)}
$$
commutes, where $\alpha_V$ and $\alpha$ are the source projections. Now
$T(\alpha_V)(Z,X) = X$ and so it follows that $T(b)(X) = 0$. Similarly,
$Z\in T^bM$. It remains to prove that $A(\Tilde{c})(Z) = A(\Tilde{c})(X)$.

Taking $S = R(\Tilde{b})$ with its horizontal structure over $R(c)$, the
source and target maps $\tilbeta_V\colon S\to R(b),\ (x,y,z,w)\mapsto (w,z),$
and $\tilalpha_V\colon (x,y,z,w)\mapsto (y,x)$, are morphisms. Applying the
Lie functor, we have $A(\tilalpha_V), A(\tilbeta_V)\colon
T^{\Tilde{b}}(R(c))\to T^bM$ and $A(\tilalpha_V)(Z,X) = X,\
A(\tilbeta_V)(Z,X) = Z.$ Now all that is necessary is to observe that on
the groupoid level we have
$\Tilde{c}\circ \tilalpha_V = \Tilde{c}\circ \tilbeta_V$.

This proves that $T^{\Tilde{b}}(R(c))\subseteq R(A(\Tilde{c})).$ Since the
maps $c,b,\Tilde{c},\Tilde{b}$ are all submersions, a dimension count
shows that equality holds.
\boom

This is of course the calculation of $A_HS$ as in
\cite[\S4, p.222]{Mackenzie:1992}. In this case however the Lie algebroid
structure on base $R(c)$ is that of an involutive distribution on $R(c)$
(as well as an action Lie algebroid), and the groupoid structure on base
$T^bM$ is the kernel pair of $A(\Tilde{c})$ (as well as an action groupoid).
Now we have
$$
A^2 = T^{A(\Tilde{c})}T^bM.
$$
This is an involutive distribution on $T^bM$ but we prefer to regard it as a
submanifold of the double vector bundle $T^2M$. It is easy to see from the
construction, or directly, that $T(p)$ maps $A^2$ to $T^cM$, and we
consequently have a double vector bundle
\begin{equation}                    \label{dvb:A2}
\dvb{A^2}{T^bM}{T^cM}{M.}{T(p)}{p_{T}}{p}{p}
\end{equation}

\begin{thm}                              \label{thm:vac}
The map $(T(p),p_T)\colon A^2\to T^cM\oplus T^bM$ is a diffeomorphism.
\end{thm}

The proof of the following lemma is straightforward.

\begin{lem}
If $\phi\colon \Om'\to \Om$ is a morphism of Lie groupoids over
$f\colon M'\to M$, and $\phi^!\colon\Om'\to f^!\Om$ is a diffeomorphism,
then $A(\phi)^!\colon A\Om\to f^!A\Om$ is also a
diffeomorphism.
\end{lem}

\noindent{\sc Proof of \ref{thm:vac}:}\
We know that $R(\Tilde{b})\to R(c)*R(b),\ (x,y,z,w)\mapsto((z,x),(y,x)),$
is a diffeomorphism, so $\tilalpha_V\colon R(\Tilde{b})\to R(b),\
(x,y,z,w)\mapsto (y,x),$ which is a morphism over
$R(c)\to M,\ (z,x)\mapsto x,$ satisfies the condition of the lemma. Hence
$$
T^{\Tilde{b}}(R(c))\to R(c)*T^bM,\qquad (Z_z,X_x)\mapsto ((z,x),X),
$$
is a diffeomorphism. Changing the point of view, it follows that the
groupoid morphism $R(A(\Tilde{c}))\to R(c),\ (Z_z,X_x)\mapsto (z,x),$
over $T^bM\to M$ satisfies the condition of the lemma. Hence
$A^2\to T^bM\oplus T^cM$ is a diffeomorphism.
\boom

Before proceeding, we recall some very basic facts about connections
\cite[XVII\S18]{Dieudonne:III}. Associated with the double vector bundle
$T^2M$ are the two exact sequences
$$
p^!TM \inj T^2M {\buildrel T(p)^! \over \surj} p^! TM,\qquad
p^!TM \inj  T^2M {\buildrel p_T^!\over \surj} p^! TM
$$
where the central terms are respectively the vector bundles
$p_T\colon T^2M\to TM$ and $T(p)\colon T^2M\to TM$. A connection in $M$
is a map ${\sf C}\colon TM\oplus TM\to T^2M$ which is simultaneously a linear
right--inverse for both sequences. We take it that $T(p)({\sf C}(X,Y)) = Y$
and $p_T({\sf C}(X,Y)) = X$. Given a connection ${\sf C}$ and a vector field
$X$ on $M$, define a vector field $X^{\sf C}$ on $TM$ by $X^{\sf C}(Y) =
{\sf C}(Y,X)$. The connection is {\em flat} if $[X,Y]^{\sf C} =
[X^{\sf C},Y^{\sf C}]$ for all $X,Y\in{\cal X}(M)$.

Given $X,Y\in{\cal X}(M)$ and $m\in M$, consider
$T(Y)(X(m)) - X^{\sf C}(Y(m))$. This is a vertical tangent vector at $Y(m)$
and so corresponds to an element of $T_m(M)$, which is denoted
$\nabla_X(Y)(m)$. This defines the associated {\em Koszul connection}
$\nabla$. There is a bijective correspondence between connections and Koszul
connections.

Given a connection ${\sf C}$, the {\em conjugate connection} ${\sf C}'$
is ${\sf C}' = J\circ {\sf C}\circ J_0$ where $J$ is the canonical involution
in $T^2M$ and $J_0\colon TM\oplus TM\to TM\oplus TM$ interchanges the
arguments. The corresponding Koszul connection $\nabla'$ is given by
$\nabla'_X(Y) = \nabla_Y(X) + [X,Y]$.

We now consider the inverse of the diffeomorphism in Theorem~\ref{thm:vac}
as constituting a \lq\lq bipartial\rq\rq\ connection in $M$ adapted to the
two foliations $R(c)$ and $R(b)$. More precisely, it is a partial connection
in the vector bundle $T^bM$ adapted to $R(c)$. The above observations about
connections in $M$ apply with the obvious modifications.

Let ${\cal X}^c$ and ${\cal X}^b$ denote the modules of sections of $T^cM$
and $T^bM$. Given $X\in{\cal X}^c$ there is a unique
$\Bar{X}\in\Gamma_{T^bM}(A^2)$ which projects to $X$, and this induces as
above an operator $\nabla^b_X\colon {\cal X}^b\to{\cal X}^b$. We thus obtain
a \lq\lq bipartial Koszul connection\rq\rq\
$$
\nabla^b\colon {\cal X}^c\times{\cal X}^b\to{\cal X}^b.
$$
Since $T(p)\colon A^2\to T^aM$ is a Lie algebroid morphism,
we have $\Bar{[X_1,X_2]} = [\Bar{X_1},\Bar{X_2}]$ and so $A^2$ defines a flat
bipartial connection. In terms of $\nabla^b$ we have $\nabla^b_{[X_1,X_2]}
= [\nabla^b_{X_1},\nabla^b_{X_2}]$ for all $X_1,X_2\in{\cal X}^c$.

Now interchanging $c$ and $b$,
we obtain $A_2 = T^{A(\Tilde{b})}T^cM$ which is a double vector bundle
$$
\dvb{A_2}{T^bM}{T^cM}{M,}{p_T}{T(p)}{p}{p}
$$
and a sub double vector bundle of $T^2M$. From Theorem \ref{thm:jS1} we now
have:

\begin{prop}                                 \label{prop:A22}
The canonical involution $J\colon T^2M\to T^2M$ carries $A^2$
isomorphically onto $A_2$.
\end{prop}

$A_2$ can be considered to be a flat partial connection in $T^cM$ adapted
to $R(b)$. In the same way as above we obtain a flat partial Koszul
connection $\nabla^c\colon{\cal X}^b\times{\cal X}^c\to{\cal X}^c.$ The
following definition now seems reasonable.

\begin{df}
Let $M$ be a manifold and let $c\colon M\to Q_H$ and $b\colon M\to Q_V$ be
two surjective submersions. Then an {\em infinitesimal affinoid structure}
on $(M,c,b)$ is a sub double vector bundle $A^2$ of $T^2M$ of the form
(\ref{dvb:A2}) such that $(T(p),p_T)\colon A^2\to T^cM\oplus T^bM$ is a
diffeomorphism and such that $T(p)\colon A^2\to T^cM$ and
$T(p)\colon A_2 = J(A^2)\to T^bM$ are Lie algebroid morphisms.
\end{df}

Evidently any infinitesimal affinoid structure induces two partial Koszul
connections $\nabla^b$ and $\nabla^c$ as above. By Proposition
\ref{prop:A22} and the definition of conjugate connections, we have
(compare Remark~3.2 in \cite{Weinstein:1990}):

\begin{prop}                                        \label{prop:conj}
In any infinitesimal affinoid structure, $\nabla^b$ and $\nabla^c$ are
conjugate connections; that is, for all $X\in{\cal X}^c,\ Y\in{\cal X}^b$,
$$
\nabla^c_Y(X) = \nabla^b_X(Y) + [X,Y].
$$
\end{prop}

\begin{prop}                             \label{prop:matchconn}
For $X,X_1,X_2\in{\cal X}^c,\ Y,Y_1,Y_2\in{\cal X}^b$,
\begin{eqnarray*}
\nabla^b_X[Y_1,Y_2] & = & [\nabla^b_X(Y_1),Y_2] + [Y_1,\nabla^b_X(Y_2)]
                                    + \nabla^b_{\nabla^c_{Y_2}(X)}(Y_1)
                                   - \nabla^b_{\nabla^c_{Y_1}(X)}(Y_2),\\
\nabla^c_Y[X_1,X_2] & = & [\nabla^c_Y(X_1),X_2] + [X_1,\nabla^c_Y(X_2)]
                                    + \nabla^c_{\nabla^b_{X_2}(Y)}(X_1)
                                    - \nabla^c_{\nabla^b_{X_1}(Y)}(X_2).
\end{eqnarray*}
\end{prop}

\pf
Since this is a purely formal calculation, we may as well consider an ordinary
connection $\nabla$ in $M$ with conjugate connection $\nabla'$. We calculate
$R'$, the curvature of $\nabla'$.

By definition, $R'(X,Y)(Z) = \nabla'_{[X,Y]}(Z) - [\nabla'_X,\nabla'_Y](Z).$
Substituting in $\nabla'_X(Y) = \nabla_Y(X) + [X,Y]$, we obtain
\begin{eqnarray*}
R'(X,Y)(Z) & = & \nabla_Z[X,Y] - [X,\nabla_Z(Y)] - [\nabla_Z(X),Y] \\
           &   & - \nabla_{[Y,Z]}(X) - \nabla_{[Z,X]}(Y)
                 - \nabla_{\nabla_{Z}(Y)}(X) + \nabla_{\nabla_{Z}(X)}(Y).
\end{eqnarray*}
So if $R' = 0$ we get
$$
\nabla_Z[X,Y] = [X,\nabla_Z(Y)] + [\nabla_Z(X),Y]
                  + \nabla_{\nabla'_Y(Z)}(X) - \nabla_{\nabla'_X(Z)}(Y).
$$

In the case of the bipartial flat connections of an infinitesimal affinoid we
can apply this calculation to both $\nabla^b$ and $\nabla^c$.
\boom

The equations of \ref{prop:conj} and \ref{prop:matchconn} show that $T^cM$
and $T^bM$ are a matched pair of Lie algebroids in the sense of Mokri
\cite{Mokri:1997}. The corresponding representations are $\nabla^c$ and
$\nabla^b$. Indeed $(T^cM, T^bM)$ is the matched pair of Lie algebroids
arising from the vacant double groupoid structure on $S$.

It is now clear that we could equivalently define an infinitesimal affinoid
structure on $(M,c,b)$ to be a pair of partial flat connections $\nabla^c$
and $\nabla^b$ which satisfy the equations of Proposition
\ref{prop:matchconn}.

We now turn to the infinitesimal form of butterfly diagrams, which is easily
defined.

\begin{df}
An {\em infinitesimal butterfly diagram} is a diagram of Lie algebroid
morphisms of the form shown in Figure~\ref{fig:ibd},
\begin{figure}[htb]
\begin{picture}(200,180)(-100,0)
\put(20,150){$R_b$}                   \put(155,150){$R_c$}
\put(95,85){$L$}
\put(20,20){$A_h$}                     \put(160,20){$A_v$}
\put(45,140){\vector(1,-1){40}}
\put(25,140){\vector(0,-1){100}}
\put(115,75){\vector(1,-1){40}}
\put(165,140){\vector(0,-1){100}}
\put(150,140){\vector(-1,-1){40}}
\put(80,75){\vector(-1,-1){40}}
\put(10,90){$\Tilde{c}$}
\put(180,90){$\Tilde{b}$}
\put(50,110){$i_v$}
\put(50,60){$\Bar{c}$}
\put(140,110){$i_h$}
\put(140,60){$\Bar{b}$}
\end{picture}
\caption{\ \label{fig:ibd}}
\end{figure}
\noindent such that the two vertical arrows are action morphisms over
surjective submersions $c\colon M\to Q_H$ and $b\colon M\to Q_V$, the two
upper diagonals are base-preserving embeddings over $c$ and $b$, the two
lower diagonals are inductors, and the two sequences are exact at their
central term.
\end{df}

Applying the Lie functor to a butterfly diagram clearly leads to an
infinitesimal butterfly diagram, since the Lie functor sends action
morphisms to action morphisms, inductors to inductors, embeddings to
embeddings, and is exact \cite{HigginsM:1990a}.

We now show that an infinitesimal affinoid structure gives rise to an
infinitesimal butterfly diagram, provided that the fibres of
$c$ and $b$ are simply--connected. The first step holds without any
connectivity assumption. Its proof is a straightforward calculation.
(A more general result of this type is given in \cite{Mokri:1997}.)

\begin{prop}                                \label{prop:diag}
Let $(M,c,b)$ have an infinitesimal affinoid structure with associated
partial connections $\nabla^c$ and $\nabla^b$. Denote by $L$ the vector
bundle direct sum $T^cM\oplus T^bM$. Then $L$ has a Lie algebroid
structure over $M$ with anchor $a\colon L\to TM$ given by
$a(X\oplus Y) = X + Y$, and bracket
$$
[X\oplus Y,X'\oplus Y'] =
             \{[X,X'] + \nabla^c_Y(X') - \nabla^c_{Y'}(X)\}\oplus
             \{[Y,Y'] + \nabla^b_X(Y') - \nabla^b_{X'}(Y)\}
$$
for $X,X'\in{\cal X}^c,\ Y,Y'\in{\cal X}^b$.
\end{prop}

Until Theorem~\ref{thm:iasibd}, consider an infinitesimal affinoid structure
on $(M,c,b)$ for which the fibres of $c$ and $b$ are simply connected. The
partial connections $\nabla^c$ and $\nabla^b$ can be considered as Lie
algebroid morphisms $T^bM\to\CDO(T^cM)$ and $T^cM\to\CDO(T^bM)$ and
accordingly integrate to give linear actions of $R(b)$ on the vector bundle
$T^cM$ and of $R(c)$ on the vector bundle $T^bM$. Denote these actions by
$\theta^c$ and $\theta^b$.

We need to show that $L$ quotients over $T^bM$ and $T^cM$ to give Lie
algebroids $A_v$ on $Q_V$ and $A_h$ on $Q_H$. We follow the method of
\cite[\S4]{HigginsM:1990a}. Say that $X\in{\cal X}^c$ is
$\theta^c$-{\em stable} if
$\theta^c(y,x)X(x) = X(y)$ for all $(y,x)\in R(b)$. This is equivalent to
the condition that $X$ be $\nabla^c$-parallel, that is, that
$\nabla^c_Y(X) = 0$ for all $Y\in{\cal X}^b$. We also say that
$X\oplus Y\in \Gamma L$ is $\theta^c$-stable or $\nabla^c$-parallel if
$X$ is so.

To show that $T^bM$ is a $b$-ideal of $L$ we need to verify the following
three conditions:
\begin{enumerate}
\item If $X\oplus Y,\ X'\oplus Y'$ are $\theta^c$-stable, then
$[X\oplus Y, X'\oplus Y']$ is also;
\item If $X\oplus Y\in\Gamma L$ is $\theta^c$-stable, and $Y'\in{\cal X}^b$,
then $[X\oplus Y, 0\oplus Y']$ is in $T^bM$;
\item The map $L/T^bM\to TM/T^bM$ induced by the anchor of $L$ is equivariant
with respect to $\theta^c$ and the natural action of $R(b)$ on
$TM/T^bM\isom b^!TQ_V.$
\end{enumerate}

The first two conditions are easily checked. For the third, note that the
natural action of $R(b)$ on $TM/T^bM$ differentiates to $D_Y(\Bar{Z}) =
\Bar{[Y,Z]}$ where $Y\in{\cal X}^b,\ Z\in{\cal X}(M)$; here the bar
denotes the class modulo $T^bM$. It therefore suffices to check that
$\Bar{\nabla^c_Y(X)} = \Bar{[Y,X]}$ and this follows from the fact that
$\nabla^c_Y(X) - [Y,X] = \nabla^b_X(Y)\in{\cal X}^b.$

From \cite[4.5]{HigginsM:1990a} it therefore follows that the vector bundle
$T^cM\isom L/T^bM$ descends to a vector bundle $A_v$ on $Q_V$; that is,
$T^cM$ is the vector bundle pullback of $A_v$ over $b$. Denote the map
$T^cM\to A_v$ by $\Tilde{b}$ and the map $L\to L/T^bM\to A_v$ by $\Bar{b}$.
Further, $A_v$ has a Lie algebroid structure over $Q_V$ with respect to which
$\Tilde{b}$ (and therefore $\Bar{b}$) is a Lie algebroid morphism. The
bracket of two sections $x,x'\in\Gamma A_v$ is obtained by taking the inverse
image sections $X,X'\in\Gamma L$; by (i) above, $[X,X']$ is $\theta^c$-stable
and therefore descends to a section $[x,x']$ of $A_v$.

The rank of $A_v$ is the same as that of $T^cM$ and so $\Bar{b}$ is an action
morphism. Since $\Tilde{b}$ has kernel $T^bM$ by construction, it is an
inductor. Carrying out the same construction with $T^cM$ as kernel, we
have proved the following.

\begin{thm}                                      \label{thm:iasibd}
Let $(M,c,b)$ have an infinitesimal affinoid structure and assume that $c$
and $b$ have simply connected fibres. Then the above construction yields an
infinitesimal butterfly diagram.
\end{thm}

Conversely consider an infinitesimal butterfly diagram as in
Figure~\ref{fig:ibd}. An inductor of Lie algebroids is essentially a
pullback in the category of Lie algebroids, and so its kernel is the
vertical bundle of the base map. Thus we have $R_b = T^bM$ and
$R_c = T^cM$.

Since $\Bar{c}\colon L\to A_h$ is an inductor over $c$, there is an
isomorphism of Lie algebroids
$$
L\to TM\pback{c^!TQ_H}c^!A_h,\qquad Z\mapsto (a_L(Z),\Bar{c}^!(Z)),
$$
where $a_L$ is the anchor of $L$. And since $\Tilde{c}\colon T^bM\to A_h$
is an action morphism, we know that $\Tilde{c}^!\colon T^bM\to c^!A_h$
is an isomorphism; denote its inverse by $\eta$. Now it is readily checked
that
$$
L\to T^cM\oplus T^bM,\qquad
   Z\mapsto (a_L(Z) - \eta(\Bar{c}^!(Z)))\oplus \eta(\Bar{c}^!(Z))
$$
is an isomorphism of vector bundles, and that $i_h$ and $i_v$ are now
represented by $X\mapsto X\oplus 0$ and $Y\mapsto 0\oplus Y$. We can
therefore apply the following lemma, whose proof is purely formal.

\begin{lem}                                   \label{lem:sum}
Let the vector bundle direct sum $T^cM\oplus T^bM$ have a Lie algebroid
structure over $M$ with respect to which $T^cM$ and $T^bM$ are Lie
subalgebroids. Define
$\nabla^c\colon{\cal X}^b\times{\cal X}^c\to {\cal X}^c$ and
$\nabla^b\colon{\cal X}^c\times{\cal X}^b\to {\cal X}^b$ by
$$
[0\oplus Y,X'\oplus 0] = \nabla^c_Y(X')\oplus -\nabla^b_{X'}(Y).
$$
Then $\nabla^c$ and $\nabla^b$ are flat partial connections and satisfy
the relations in Proposition \ref{prop:matchconn}.
\end{lem}

This completes the proof of the following result.

\begin{thm}                                      \label{thm:ibdias}
Let $L$ be an infinitesimal butterfly diagram over $c\colon M\to Q_H$ and
$b\colon M\to Q_V$. Then the above construction yields an infinitesimal
affinoid structure on $M$.
\end{thm}

Theorems \ref{thm:iasibd} and \ref{thm:ibdias} establish an equivalence
between infinitesimal affinoid structures and infinitesimal butterfly
diagrams when the fibres of $c$ and $b$ are simply connected.

Infinitesimal butterfly diagrams behave in some respects like Morita
equivalences: for example, assuming that the fibres of $c$ and $b$ are
simply--connected, one can easily prove that $A_h$ is integrable if and
only if $A_v$ is integrable. However, given an arbitrary Lie algebroid $A$,
there appears to be no canonical construction of an infinitesimal butterfly
diagram with $A_h = A_v = A$, and this indicates severe limitations for the
concept.

Lastly in this section consider a generalized principal bundle
$P(B,G,p)(Q_V,f)$ \cite{Kock:1988} as defined in \cite[3.3]{Mackenzie:1992}.
Thus $G\gpd Q_V$
is a Lie groupoid acting smoothly and freely to the right on a surjective
submersion $f\colon P\to Q_V$ with quotient manifold $B = P/G$ and
projection $p\colon P\to B$. This defines an affinoid structure on
$P$ with respect to $c = p\colon P\to B$ and $b = f\colon P\to Q_V$.
The vertical groupoid $G_v$ identifies canonically with the original $G$.
For clarity, denote the elements of $G_h$ by $\langle y,x\rangle$ where
$x,y\in P,\ f(x) = f(y),$ and $\langle y,x\rangle = \langle yg,xg\rangle$
for any $g\in G$ with $\beta g = f(x)$.

Since $P(B,G,p)(Q_V,f)$ is presented as a generalization of the
concept of principal bundle, it is reasonable to extend the notion of
Atiyah sequence to it. (The omitted details in what follows may be found
by extending the account in \cite[App.A]{Mackenzie:LGLADG}.) Consider the
vertical tangent bundle $T^fP$. The action of $G$ on $P$ lifts to a
right action of $G$ on $T^fP$ and remains free. Denote elements of
$T^fP/G$ by $\langle X_x\rangle$ where $X_x\in T^f_xP$. Then
$T^fP/G$ is a vector bundle over $B$; if $\langle X_x\rangle$ and
$\langle Y_y\rangle$ have $p(x) = p(y)$, then there exists $g\in G$ with
$y = xg$ and we define $\langle X_x\rangle + \langle Y_y\rangle
= \langle Xg + Y\rangle$, as in the standard case. An $f$-vertical
vector field $X$ on $P$ may be defined to be {\em $G$--invariant} if
$X(xg) = X(x)g$ for all $x\in P$ and $g\in G$ with $f(x) = \beta g$;
one obtains a $C(B)$--module of $G$--invariant vector fields which is
in bijective correspondence with the module of sections of $T^fP/G\to B$.
Now the bracket of $G$--invariant vector fields transfers to
$\Ga(T^fP/G)$ and makes it a Lie algebroid with anchor the quotient to
$T^fP/G\to TB$ of $T(p)\colon T^fP\to TB$. This might be called the
{\em generalized Atiyah sequence of} $P(B,G,p)(Q_V,f)$. As in the standard
case, $T^fP\to T^fP/G$ is an action morphism of Lie algebroids over $p$.

\begin{prop}                                     \label{prop:gAs}
The generalized Atiyah sequence just constructed is canonically
isomorphic to $A(G_h)$.
\end{prop}

\pf
The morphism $\Tilde{c}\colon R(f)\to G_h$ over $p\colon P\to B$ is
$(y,x)\mapsto\langle y,x\rangle.$ It induces a Lie algebroid morphism
$A(\Tilde{c})\colon T^fP\to A(G_h)$ which is constant on the orbits of
$G$ and therefore induces a Lie algebroid morphism $R\colon T^fP/G\to A(G_h)$
over $B$. Since $\Tilde{c}$ is an action morphism, it is a fibrewise
diffeomorphism, and this property is inherited by $A(\Tilde{c})$ and $R$.
Since $R$ is base--preserving, it is therefore an isomorphism of Lie
algebroids.
\boom

Note that the method available in the standard case
\cite[III~3.20]{Mackenzie:LGLADG} cannot be used here, since $P$ cannot
generally be embedded in $G_h$.

Proposition \ref{prop:gAs} rules out much interest in the notion of
generalized Atiyah sequence. Since any Lie groupoid $G$ acts freely to the
right on its target projection, yielding a generalized principal bundle
whose $G_h$ is again canonically isomorphic to $G$, any Lie algebroid which
is the Lie algebroid of a Lie groupoid may be constructed as a generalized
Atiyah sequence. There is no prospect of extending to generalized
Atiyah sequences the very rich theory known for the usual notion of
Atiyah sequence, which depends essentially on its transitivity.

Although generalized principal bundles are equivalent to both affinoid
structures and butterfly diagrams, it is clear that generalized Atiyah
sequences embody only a part of the information in an infinitesimal
butterfly diagram.

\section{DUALITY AND THE COTANGENT DOUBLE \break GROUPOID}
\label{sect:datcdg}

As a preliminary to the next section, we show that the cotangent of a
double Lie groupoid itself has a double groupoid structure. We begin
by recalling the notion of dual for \VBgpd s due to Pradines
\cite{Pradines:1988}.

Consider a \VBgpd\ $(\Om;G,A;M)$ as in \cite[\S4]{Mackenzie:1992}, with core
$K = \{\xi\in\Om\st \tilalpha(\xi) = 0_m,\ \tilq(\xi) = 1_m,\ \exists m\in M\}$.
The vector bundle operations on $\Om$ restrict to give $K$ the structure of a
vector bundle over $M$.

Let $\Om^*$ be the dual of $\Om$ as a vector bundle over $G$. Define a
groupoid structure on $\Om^*$ with base $K^*$ as follows. Take $\Phi\in
\Om^*_g$ where $g\in G_m^n$. Then the source and target of $\Phi$ in
$K^*_m$ and $K^*_n$ respectively are
\begin{equation}                                    \label{eq:dual1}
\langle\tilalpha_*(\Phi),\kappa\rangle =
         \langle\Phi,-\til0_g\kappa^{-1}\rangle, \quad \kappa\in K_m,
\qquad\qquad
\langle\tilbeta_*(\Phi),\kappa\rangle =
         \langle\Phi, \kappa\til0_g\rangle, \quad \kappa\in K_n.
\end{equation}
For the composition, take $\Psi\in\Om^*_h$ with
$\tilalpha_*(\Psi) = \tilbeta_*(\Phi)$. Any element of $\Om_{hg}$ can
be written as a product $\eta\xi$ where $\eta\in\Om_h$ and $\xi\in\Om_g$.
Now the compatibility condition on $\Psi$ and $\Phi$ ensures that
\begin{equation}                                    \label{eq:dual2}
\langle\Psi\Phi, \eta\xi\rangle = \langle\Psi,\eta\rangle +
                                    \langle\Psi,\xi\rangle
\end{equation}
is well defined. The identity element of $\Om^*$ at $\theta\in K^*_m$ is
$\tilone_\theta\in\Om^*_{1_m}$ defined by
\begin{equation}                                    \label{eq:dual3}
\langle\tilone_\theta,\tilone_X + \kappa\rangle = \langle\theta,\kappa\rangle,
\end{equation}
where any element of $\Om_{1_m}$ can be written as $\tilone_X + \kappa$
for some $X\in A_m$ and $\kappa\in K_m$.

It is straightforward to check that $(\Om^*;G,K^*;M)$ is a \VBgpd, the
{\em dual \VBgpd\ to $\Om$}. The core of $\Om^*$ is the vector bundle
$A^*\to M$, with the core element corresponding to $\phi\in A^*_m$ being
$\Bar\phi\in\Om^*_{1_m}$ defined by
\begin{equation}
\langle\Bar\phi,\tilone_X + \kappa\rangle =
   \langle\phi,X + \da_A(\kappa)\rangle
\end{equation}
for $X\in A_m,\ \kappa\in K_m$. Here $\da_A\colon K\to A$ is the map defined
by $\tilbeta\colon\Om\to A$ as in \cite[\S5]{Mackenzie:1992}.

Now consider a morphism of \VBgpd s which preserves the lower groupoids:
$$
(F;\id_G,f;\id_M)\colon(\Om;G,A,M)\to (\Om';G,A';M)
$$
and denote the core morphism $K\to K'$ by $f_K$. The proof of the following
result is simple.

\begin{prop}                                  \label{prop:dualmorph}
The dual morphism $F^*\colon\Om'^*\to\Om^*$ is a morphism of the dual
\VBgpd s, with base map $f_K^*\colon K'^*\to K^*$ and core morphism
$f^*\colon A'^*\to A^*$.
\end{prop}

For any Lie groupoid $G\gpd M$, the dual of the tangent \VBgpd\
$(TG;G,TM;M)$ is the cotangent groupoid $(T^*G;G,A^*G;M)$; the conventions
used above give precisely the structure given in \cite[\S7]{MackenzieX:1994}.
If now $G\gpd M$ is a Poisson groupoid, $\pi_G^\#\colon T^*G\to TG$ is a
morphism of \VBgpd s with respect to the base maps $\id_G$ and
$a_*\colon A^*G\to TM$ (see \cite{AlbertD:1991pp} or
\cite[8.1]{MackenzieX:1994}). It now follows from \ref{prop:dualmorph} and
from the skewsymmetry of $\pi_G^\#$, that the core map of $\pi_G^\#$ is
$-a_*^*$.

The duality for \VBgpd s can be extended to a duality between \LAgpd s and
what one might call ${\cal PVB}$--groupoids: \VBgpd s $(\Om;G,E;M)$ in which
$\Om$ and $E$ carry Poisson structures making $\Om\gpd E$ a Poisson groupoid.
This notion will be developed elsewhere.

We now give a double version of the cotangent groupoid. Consider a double Lie
groupoid $(S;H,V;M)$ with core groupoid $C\gpd M$ as in \cite{BrownM:1992}
or \cite[\S2, p.197]{Mackenzie:1992}. Denote by $A_HS$ and $A_VS$ the Lie
algebroids of the two groupoid structures on $S$.

The standard cotangent groupoid structure (\cite{CDW} or \cite{Pradines:1988})
gives $T^*S$ groupoid structures over both $A^*_HS$ and $A^*_VS$. In turn,
since both $A_HS$ and $A_VS$ are \LAgpd s with core $AC$, as in
\ref{prop:cores}, the duals $A^*_HS$
and $A^*_VS$ both have groupoid structures on base $A^*C$.

\begin{thm}
With the structures just described, $T^*S$ is a double Lie groupoid
\begin{equation}                       \label{eq:doublecot}
\matrix{&&{\tilalpha_{*H},\tilbeta_{*H}}&&\cr
          &T^*S&\lgpd &A^*_HS&\cr
          &&&&\cr
          {\tilalpha_{*V},\tilbeta_{*V}}&\vgpd&&\vgpd&{\alpha_{*V},\beta_{*V}}\cr
          &&&&\cr
          &A^*_VS&\lgpd &A^*C&\cr
          &&{\alpha_{*H},\beta_{*H}}&&\cr}
\end{equation}
\end{thm}

\pf
Using the groupoid structures defined in (\ref{eq:dual1})---(\ref{eq:dual3})
on the side groupoids and, on $T^*S$, the ordinary cotangent groupoid
structures which are a special case of these, this is a long but
straightforward verification.
\boom

\begin{thm}
The core groupoid of (\ref{eq:doublecot}) is naturally isomorphic with the
cotangent groupoid $T^*C\gpd A^*C$ of the core of $S$.
\end{thm}

\pf
Define a map $E\colon T_CS\to TC$ which sends each vector tangent to $S$
at a point of the core, to a vector at the same point tangent to the core
itself. Take $c\in C$ with $v = \da_V(c),\ h = \da_H(c)$, where $v\in V_m,\
h\in H_m$. Represent an element $\xi\in T_cS$ as
$$
\matrix{w&\raise.7ex\hbox{$W$}&x\cr
                  \raise.35in\hbox{$Z$}&\ssq&\raise.35in\hbox{$X$}\cr
                  z&Y&y\cr}\hskip-.7in\raise5pt\hbox{$\xi$}
$$
where $X = T(\tilalpha_H)(\xi)\in T_{1^V_m}V,\ x = T(\beta_V)(X) =
T(\alpha_H)(W) \in T_mM$, et cetera, and define
\begin{equation}
E(\xi) = \xi - T(L^V_c)T(\tilone^H)(X - T(1^V)(x))
             - T(L^H_c)T(\tilone^V)(Y - T(1^H)(z)),
\end{equation}
where $L^V_c$ and $L^H_c$ denote left translation in the two groupoid
structures on $S$.

Take $\sigma\in T^*_cC$ and define $\Sigma\in T^*_cS$ by
$\Sigma(\xi) = \sigma(E(\xi)),\ \xi\in T_cS$. We must show that $\Sigma$
is a core element. Denote $\tilalpha_{*V}(\Sigma)\in A^*_VS|_{1^H_m}$ by
$\rho$, and $\alpha_{*H}(\rho)\in A^*_mC$ by $\theta$. Then for all
$\kappa\in A_mC$,
$$
\langle\theta,\kappa\rangle = -\rho(T(\tiliota_H)(\kappa))
= \Sigma(T(L^V_c)T(\tiliota_V)T(\tiliota_H)(\kappa))
$$
where $\tiliota_H$ and $\tiliota_V$ are the two inversions in $S$.

We first prove that $\rho = 1^{*V}_\theta$. Take any
$\xi = A(\tilone^H)(X) + \kappa\in A_VS|_{1^H_m}$ where $X\in A_mV,\
\kappa\in A_mC$. Then $1^{*V}_\theta(\xi) = \langle\theta,\kappa\rangle$
and $\rho(\xi) = - \Sigma(T(L^V_c)T(\tiliota_V)(\xi))$. Now although it is
natural to consider $\kappa$ as an element of $A_VS$ in this context, it may
also be regarded as an element of $A_HS$ and we then get
$$
T(\tiliota_H)(\kappa) = T(\tilone^H)(Z) - \kappa
$$
where $Z = \da_{AV}(\kappa)$. Next, applying $E$ to
$T(L^V_c)T(\tiliota_V)T(\tilone^H)(Z)$, or equally with $X$ in place of $Z$,
gives $0\in T_cC$, where we use
$\tiliota_V\circ\tilone^H = \tilone^H\circ\iota_V$. Putting these together,
the equation for $\rho$ follows.

We must also prove that $\theta\in A^*_mC$ is the source of $\sigma$ with
respect to $T^*C\gpd A^*C$. For $\kappa\in A_mC$ the above shows that
$$
\langle\theta,\kappa\rangle = -\Sigma(T(L^V_c)T(\tiliota_V)(\kappa)).
$$
Calculating $E$ of the argument on the right hand side, we obtain
\begin{equation}                             \label{eq:sum}
T(L^V_c)T(\tiliota_V)(\kappa) + T(L^H_c)T(\tilone^V)(T(1^H)(a_HW) - W),
\end{equation}
where $W = \da_{AH}(\kappa)$. Write $\kappa = \ddt{c_t}$ where the curve
$c_t$ in $C$ is of the form
$$
\matrix{ &\raise.7ex\hbox{$h_t$}&\cr
                  \raise.35in\hbox{$v_t$}&\ssq&\raise.35in\hbox{$1^V_m$}\cr
                   &1^H_m&\cr}\hskip-.7in\raise5pt\hbox{$c_t$}
$$
Then the first term in (\ref{eq:sum}) is $\ddt{c\,\vcomp\, c_t^{-V}}$ and the
second is $\ddt{c\,\hcomp\, \tilone^V_{h_t^{-1}}}$. The sum is therefore the
derivative of the product shown in Figure~\ref{fig:squares},
\begin{figure}[htb]
\begin{picture}(350,120)  
\put(160,65){$\ssq$}                    \put(230,65){$\ssq$}
\put(180,85){$c$}                       \put(255,85){$\tilone^V_{h_t^{-1}}$}
\put(160,0){$\ssq$}                     \put(230,0){$\ssq$}
\put(180,25){$c_t^{-V}$}                \put(255,25){$\tilone^V_{h_t^{-1}}$}
\end{picture}\caption{\ \label{fig:squares}}
\end{figure}
and in terms of the multiplication $\dcomp$ and the inversion
$c\mapsto c^{-C}$ in $C\gpd M$, this product is $c\dcomp c_t^{-C}$
(see \cite[\S2, p.197]{Mackenzie:1992} or \cite{BrownM:1992}). The
derivative is therefore $T(L^C_c)T(\iota_C)(\kappa)$ and so we have
$\langle\theta,\kappa\rangle = \langle\alpha_{*C}(\sigma),\kappa\rangle$.

The proof that $\tilalpha_{*H}(\Sigma) = 1^{*H}_\theta$ now follows in
the same way. One likewise checks that $\beta_{*C}(\sigma) =
\beta_{*H}(\tilbeta_{*V}(\Sigma))$ and that the multiplications
correspond.

This proves that $\sigma\mapsto\Sigma$ is an isomorphism into the core of
$T^*S$. We leave the reader to check that the image is the whole of $T^*C$.
\boom

In point of fact, $T^*S$ is a triple structure, each of its four spaces being
also a vector bundle over the corresponding space of $(S;H,V;M)$, and all the
groupoid structure maps being vector bundle morphisms.

\section{SYMPLECTIC DOUBLE GROUPOIDS AND POISSON DOUBLE GROUPOIDS}
\label{sect:sdgapdg}

In introducing the concept of Poisson groupoid in \cite{Weinstein:1988},
Weinstein defined two Poisson groupoids on the same base manifold to be
{\em dual} if the Lie algebroid dual of each is isomorphic to the Lie
algebroid of the other, and described a programme for showing that, at least
locally, Poisson groupoids in duality arise as the bases of a symplectic
double groupoid. In \cite{LuW:1989} he and Lu carried out this programme
globally for Poisson groups. However neither paper fully describes the
general process by which a symplectic double groupoid gives rise to a pair
of Poisson groupoids in duality.

In this section we use the apparatus of double Lie algebroids to show that
the side groupoids of a symplectic double groupoid are Poisson groupoids
whose Lie bialgebroids \cite{MackenzieX:1994} are isomorphic under a
duality, and that the core of the double groupoid provides a symplectic
realization of the double base. The double Lie algebroid structure also
provides symplectic realizations of the linearized Poisson structures.
The proof is conceptual and reasonably simple. The key is the notion of
core as developed in \cite{BrownM:1992} and \cite[\S2]{Mackenzie:1992}.

In view of the integrability of the double base of a symplectic double
groupoid, it seems likely that there exist pairs of Poisson groupoids in
duality which do not integrate globally to a symplectic double groupoid.

As with ordinary groupoids, it is actually easier to study symplectic
structures by specializing from the Poisson case. The notion of Poisson
double groupoid defined here includes (see \cite[\S3]{Mackenzie:1992})
the Poisson affinoids of \cite{DazordLSW:1991}.

In what follows we will repeatedly use the following simple result.

\begin{prop}
\begin{enumerate}
\item Let $(\phi;\phi_H,\phi_V;\phi_M)\colon (S;H,V;M)\to (S';H',V';M')$
be a morphism of double groupoids with core morphism $\phi_C\colon C\to C'$.
Then $\da_H'\circ\phi_C = \phi_H\circ\da_H$ and $\da_V'\circ\phi_C =
\phi_V\circ\da_V$.

\item Let $(\phi;\phi_G,\phi_A;\phi_M)\colon (\Om;G,A;M)\to (\Om';G',A';M')$
be a morphism of \LAgpd s with core Lie algebroid morphism
$\phi_K\colon K\to K'$. Then $\da_{AG}'\circ\phi_K = A(\phi_G)\circ\da_{AG}$
and $\da_A'\circ\phi_K = \phi_A\circ\da_A$.
\end{enumerate}
\end{prop}

\begin{df}
A {\em Poisson double groupoid} is a double Lie groupoid $(S;H,V;M)$
together with a Poisson structure on $S$ with respect to which both
groupoid structures on $S$ are Poisson groupoids.
\end{df}

For the theory of ordinary Poisson groupoids, see \cite{Weinstein:1988},
\cite{Xu:1995}, \cite{MackenzieX:1994}. For a Poisson groupoid $G\gpd P$
we take the Poisson structure on the base to be $\pi_P^\# = a_*\circ a^*
= -a\circ a_*^*.$ This convention is opposite to that used in
\cite{MackenzieX:1994}, but is necessary in order that the two core
morphisms $a^*\colon T^*P\to A^*G$ and $-a^*_*\colon T^*P\to AG$
of the cotangent \LAgpd\ be Lie algebroid morphisms: see the remark
following \ref{prop:dualmorph}.

The two Poisson groupoid structures on $S$ induce maps
$\tila_{*H}\colon A^*_VS\to TH$ and $\tila_{*V}\colon A^*_HS\to TV$, the
anchors for the Lie algebroid structures on the duals of $A_VS\to H$ and
$A_HS\to V$, with respect to which $\pi_S^\#\colon T^*S\to TS$ is a morphism
of each of the two groupoid structures on $T^*S$. From the following result
it follows that $\pi^\#_S$ is actually a morphism of double groupoids over a
map $a_{*C}\colon A^*C\to TM$.

\begin{lem}
Let $(S;H,V;M)$ and $(S';H',V';M')$ be double Lie groupoids and let
$\phi\colon S\to S',\ \phi_H\colon H\to H'$ and $\phi_V\colon V\to V'$ be
maps such that $(\phi,\phi_H)$ and $(\phi,\phi_V)$ are morphisms of the
two ordinary groupoid structures on $S$ and $S'$. Then there is a unique
map $\phi_M\colon M\to M'$ such that $(\phi;\phi_H,\phi_V;\phi_M)$ is a
morphism of double groupoids.
\end{lem}

\pf
Take $m\in M$. The double identity $1^2_m$ can be written both as
$\tilone^V_{1^H_m}$ and as $\tilone^H_{1^V_m}$. Its image under $\phi$
is therefore an identity for both top structures on $S'$, and must
therefore be a double identity $1^2_{\phi_M(m)}$. It also follows that
$\phi_H,\ \phi_V$ and $\phi_M$ commute with the source and target
projections.

Now since $\phi(\tilone^V_h) = \tilone^V_{\phi_H(h)}$ for all $h\in H$
and $\tilone^V_{h_1}\,\hcomp\,\tilone^V_{h_2} = \tilone^V_{h_1h_2}$
for all compatible $h_1,h_2\in H$, it follows that $\phi_H$ is a morphism
over $\phi_M$. Similarly for $\phi_V$.
\boom

Returning to the Poisson double groupoid $S$, the bases $H$ and $V$ acquire
Poisson structures which we take to be $\pi_H^\# = \tila_{*H}\circ\tila^*_V$
and $\pi_V^\# = \tila_{*V}\circ\tila^*_H$. We will prove below that $H$ and
$V$ are Poisson groupoids with respect to these structures.

First note that the core morphism $T^*C\to TC$ of $\pi_S^\#$ defines a
Poisson structure on $C$; denote this by $\pi_C^\#$. The core is a morphism
of groupoids over $a_{*C}\colon A^*C\to TM$ and $C\gpd M$ is therefore a
Poisson groupoid. Give $M$ the Poisson structure
$\pi_M^\# = a_{*C}\circ a^*_C$ induced from $C$.

Since $S_V\gpd H$ is a Poisson groupoid, its Lie algebroid dual $A_V^*S$ has
a Lie algebroid structure with anchor $\tila_{*H}$. Similarly $A^*C$ has a
Lie algebroid structure with anchor $a_{*C}$.

\begin{prop}
With respect to these structures $(A^*_VS;H,A^*C;M)$ is an \LAgpd.
\end{prop}

\pf
It must be proved that $\alpha_{*H}$ and $\beta_{*H}$ and the groupoid
multiplication are Lie algebroid morphisms. This proceeds as in the case of
the cotangent \LAgpd\ of a Poisson groupoid: to prove $\beta_{*H}$ a Lie
algebroid morphism it suffices (\cite[6.1]{MackenzieX:1994} or
\cite[6.6]{HigginsM:1993}) to prove that the dual map
$$
AC * H\to A_VS,\qquad (\kappa,h)\mapsto \kappa\til0_h,
$$
is Poisson with respect to the linearized structures on $AC$ and $A_VS$.
\boom

We may therefore consider both $T^*S$ and $TS$ to be triple structures
as in Figure~\ref{fig:triples}: precisely, they are double groupoid objects
in the category of Lie algebroids.
\begin{figure}[htb]
\begin{picture}(350,200)  
\put(0,150){$\matrix{&&      &\cr
                      &T^*S&\sgpd &A^*_HS\cr
                      &&&\cr
                      &\svgpd          & &\svgpd          \cr
                      &&&\cr
                      &A^*_VS&\sgpd &A^*C\cr}$}

\put(40, 160){\vector(3,-4){40}}                
\put(110, 160){\vector(3,-4){40}}               
\put(40, 100){\vector(3,-4){40}}                
\put(120, 100){\vector(3,-4){30}}               

\put(75,70){$\matrix{&&      &\cr
                     &S &\sgpd &V\cr
                     &&&\cr
                     &\svgpd          &&\svgpd          \cr
                     &&&\cr
                     &H &\sgpd & M\cr}$}

\put(100,0){(a)}


\put(200,150){$\matrix{&&      &\cr
                      &TS&\sgpd &TV\cr
                      &&&\cr
                      &\svgpd          & &\svgpd          \cr
                      &&&\cr
                      &TH &\sgpd &TM \cr}$}

\put(240, 160){\vector(3,-4){40}}                
\put(310, 160){\vector(3,-4){40}}               
\put(240, 100){\vector(3,-4){40}}                
\put(320, 100){\vector(3,-4){30}}               

\put(275,70){$\matrix{&&      &\cr
                     &S &\sgpd &V\cr
                     &&&\cr
                     &\svgpd          &&\svgpd          \cr
                     &&&\cr
                     &H &\sgpd & M\cr}$}
\put(300,0){(b)}
\end{picture}\caption{\ \label{fig:triples}}
\end{figure}
We will call such structures {\em ${\cal LA}$--double groupoids}. The three
double structures which involve $T^*S$ or $TS$ we call the {\em upper faces}
or {\em upper structures}, the other three being the {\em lower faces} or
{\em lower structures}. That $\pi_S^\#$ respects all twelve ordinary
structures follows from the hypothesis and basic results for ordinary
Poisson groupoids.

In what follows we will use the following result repeatedly. The proof is
straightforward.

\begin{prop}                                   \label{prop:repeat}
\begin{enumerate}
\item The core of each upper face in an \LAdgpd\ is a Lie groupoid or Lie
algebroid over the core of the opposite lower face, and forms an \LAgpd\
with respect to the structures of the bases of the two cores.

\item The restriction of a morphism of \LAdgpd s to the core of an upper
face is a morphism of \LAgpd s.
\end{enumerate}
\end{prop}

The map $\tila_{*H}\colon A^*_VS\to TH$ is a morphism of the \LAgpd s
which form the bottom faces in Figure~\ref{fig:triples}, the other maps
being $a_{*C}\colon A^*C\to TM$ and $\id_H$. Its core map $A^*V\to AH$ is
therefore a morphism of Lie algebroids: denote it by $D_H$.

By \cite[\S5]{Mackenzie:1992}, the anchor $\tila_V\colon A_VS\to TH$ is a
morphism of \LAgpd s over $a_V\colon AV\to TM$ and $\id_H$, with core
morphism $\da_{AH}\colon AC\to AH$. Its dual $\tila_V^*\colon T^*H\to A^*_VS$
is therefore a morphism over $\da^*_{AH}\colon A^*H\to A^*C$ and $\id_H$,
with core map $a^*_V\colon T^*M\to A^*V$. It follows that $\pi^\#_H =
\tila_{*H}\circ\tila^*_V$ is a morphism of groupoids over
$a_{*C}\circ \da^*_{AH}\colon A^*H\to TM$; since it is also a morphism of
Lie algebroids over $H$, it is a morphism of \LAgpd s, and the core morphism
is $D_H\circ a^*_V\colon T^*M\to AH$. This proves the first part of the
following theorem.

\begin{thm}
With the induced structures, $H\gpd M$ is a Poisson groupoid with
$$
a_{*H} = a_{*C}\circ\da^*_{AH} = -a_V\circ D^*_H.
$$
The induced Poisson structure on $M$ coincides with that induced by $C$.
\end{thm}

For the last statement, recall from \cite[\S5, p.230]{Mackenzie:1992} that
$a_C = a_H\circ\da_{AH}$. Hence $a_{*C}\circ a^*_C =
a_{*C}\circ\da^*_{AH}\circ a^*_H = a_{*H}\circ a^*_H.$ The second equation
for $a_{*H}$ follows by noting that the core map for $\pi^\#_H$ is
$D_H\circ a^*_V$ and the negative dual of this is equal to the base map,
by the remark following \ref{prop:dualmorph}.

The same process can be carried out with $H$ and $V$ interchanged. We now
have two Lie algebroid morphisms, $D_H\colon A^*V\to AH$ and
$D_V\colon A^*H\to AV$, defined as the cores of $\tila_{*H}$ and
$\tila_{*V}$ respectively.

\begin{thm}
$D^*_V = - D_H.$
\end{thm}

\pf
We can regard $\tila_{*V}$ as the base map for the cotangent \LAgpd s
which form the top faces of Figure~\ref{fig:triples}. From the remark
following \ref{prop:dualmorph} it follows that the core map for these faces
is $-\tila^*_{*V}\colon T^*V\to A_HS$. In turn, as in \ref{prop:repeat},
the cores of the top faces of Figure~\ref{fig:triples}(a) are groupoids over
the cores of the bottom faces. The base map $A^*V\to AH$ for $-\tila^*_{*V}$
is, by \ref{prop:dualmorph}, the negative dual of the core map of
$\tila_{*V}$; that is, it is $-D^*_V$.

By the commutativity properties of the triple structures, one can see that
the base map of the core map of the top faces is the same as the core map
for the bottom faces. But the core map for the bottom faces is $D_H$.
\boom

We now specialize to symplectic double groupoids, as considered in
\cite{CDW}, \cite{Weinstein:1988} and \cite{LuW:1989}. We take the signs
in each structure to be the same.

\begin{df}
A {\em symplectic double groupoid} is a double Lie groupoid $(S;H,V;M)$
together with a symplectic structure on $S$ such that both groupoid
structures on $S$ are symplectic groupoids.
\end{df}

Since $\pi_S^\#\colon T^*S\to TS$ is an isomorphism of double groupoids,
it follows that the base maps $\tila_{*H}\colon A^*_VS\to TH$ and
$\tila_{*V}\colon A^*_HS\to TV$ are \LAgpd\ isomorphisms, and hence their
core maps $D_H\colon A^*V\to AH$ and $D_V\colon A^*H\to AV$ are Lie
algebroid isomorphisms. Further, the core map $\pi^\#_C\colon T^*C\to TC$
of $\pi^\#_S$ is an isomorphism, and so $C\gpd M$ is a symplectic groupoid.
This proves the following result.

\begin{thm}
Let $(S;H,V;M)$ be a symplectic double groupoid with core $C\gpd M$.
Then $C$ is a symplectic groupoid realizing $M$ and the Lie bialgebroids
$(AH, A^*H)$ and $(A^*V,AV)$ of the side Poisson groupoids are canonically
isomorphic.
\end{thm}

The negative dual of $\tila_{*H}$ gives an isomorphism of \LAgpd s from
$(T^*H;A^*H,H;M)$ to $A_VS$, the side maps being $D_V$ and $\id_H$. It
follows that $A_VS\gpd AV$ may be identified with the symplectic groupoid
$T^*H\gpd A^*H$, thus giving a symplectic realization of the linearized
Poisson structure on $AV$. Further, the Lie algebroid structure of
$A^2S\to AV$ is $AT^*H\isom T^*(A^*H)\isom T^*(AV)$. Similarly the vertical
structure $A^2S\to AH$ is $T^*(AH)\to AH$.

It is worthwhile reexamining Examples \ref{ex:symp} and \ref{ex:plg}
in the light of this section.

Finally, $T^*S$ is a symplectic double groupoid for any double Lie groupoid
$S$, and so we have the following result.

\begin{thm}
Let $(S;H,V;M)$ be a double Lie groupoid. Then $A^*_VS\gpd A^*C$ and
$A^*_HS\gpd A^*C$ are Poisson groupoids in duality.
\end{thm}

\newcommand{\noopsort}[1]{} \newcommand{\singleletter}[1]{#1}

\end{document}